\documentclass[lettersize,journal]{IEEEtran}
\IEEEoverridecommandlockouts
\usepackage[T1]{fontenc}
\newif\ifhavebib

\havebibtrue
\usepackage{latexsym}
\usepackage{setspace}
\usepackage{bm}
\usepackage{lipsum}
\usepackage{bbm}

\usepackage{blindtext}
\usepackage{graphicx}
\usepackage{graphicx}
\usepackage{array}
\usepackage{url}
\usepackage{algorithmic}
\usepackage[cmex10]{amsmath}
\usepackage{ifpdf}
\usepackage{amssymb,amsmath}
\usepackage{empheq}
\usepackage{stfloats}   
\usepackage{textcomp}
\usepackage{cite}
\usepackage{multirow}
\usepackage{diagbox}
\usepackage{subcaption}
\usepackage{mathrsfs}
\let\oldFootnote\footnote
\newcommand\nextToken\relax

\renewcommand\footnote[1]{%
    \oldFootnote{#1}\futurelet\nextToken\isFootnote}

\newcommand\isFootnote{%
    \ifx\footnote\nextToken\textsuperscript{,}\fi}

\usepackage[ruled,vlined,lined,ruled,boxed]{algorithm2e}
\usepackage[dvipsnames,usenames]{color}
\usepackage[normalem]{ulem}
\usepackage{amsthm}

\usepackage{graphicx}
\usepackage{ifpdf}
\usepackage{cite}
\usepackage{enumerate}
\usepackage{enumitem}

\usepackage{array}
\usepackage{mdwmath}
\usepackage{mdwtab}
\usepackage{eqparbox}
\usepackage{bm}

\usepackage{setspace}
\usepackage{booktabs}
\usepackage[subnum]{cases}
\usepackage{rotating}

\usepackage{listings} 
\definecolor{Red}{rgb}{1,0,0}
\definecolor{Blue}{rgb}{0,0,1}
\definecolor{Olive}{rgb}{0.41,0.55,0.13}
\definecolor{Green}{rgb}{0,1,0}
\definecolor{MGreen}{rgb}{0,0.8,0}
\definecolor{DGreen}{rgb}{0,0.55,0}
\definecolor{Yellow}{rgb}{1,1,0}
\definecolor{Cyan}{rgb}{0,1,1}
\definecolor{Magenta}{rgb}{1,0,1}
\definecolor{Orange}{rgb}{1,.5,0}
\definecolor{Violet}{rgb}{.5,0,.5}
\definecolor{Purple}{rgb}{.75,0,.25}
\definecolor{Brown}{rgb}{.75,.5,.25}
\definecolor{Grey}{rgb}{.5,.5,.5}

\newcommand{\boxhead}[5]{
   \pagestyle{myheadings}
   \thispagestyle{plain}
   \setcounter{page}{1}
   \noindent
   \begin{center}
   \framebox{
      \vbox{\vspace{2mm}
    \hbox to 6.28in { {\bf #1 \hfill} }
       \vspace{6mm}
       \hbox to 6.28in { {\Large \hfill \bf #2  \hfill} }
       \vspace{6mm}
       \hbox to 6.28in { {\it #3 #4 \hfill  #5} }
      \vspace{2mm}}
   }
   \end{center}
   \markboth{#5 -- #2}{#5 -- #2}
   \vspace*{4mm}
}

\usepackage[colorlinks=true, urlcolor=black,  linkcolor=black, citecolor=black]{hyperref}
\theoremstyle{definition}
\newtheorem{theorem}{Theorem}

\newtheorem{lemma}{Lemma}

\newtheorem{remark}{Remark}

\newtheorem{definition}{Definition}

\DeclarePairedDelimiterX{\infdivx}[2]{(}{)}{%
	#1\;\delimsize\|\;#2%
}

\DeclarePairedDelimiter{\norm}{\lVert}{\rVert}
\DeclarePairedDelimiter{\abs}{\lvert}{\rvert}
\DeclareMathOperator*{\argmax}{\mathop{\arg\max}}
\DeclareMathOperator*{\argmin}{\mathop{\arg\min}}
\def\tr{\mathop{\rm tr}\nolimits}%
\def\diag{\mathop{\rm diag}\nolimits}%
\def\rank{\mathop{\rm rank}\nolimits}%
\renewcommand{\Pr}{\mathscr{P}}

\newcommand{\av}{{\bf a}}
\newcommand{\bv}{{\bf b}}

\newcommand{\Cv}{{\bf C}}

\newcommand{\Dv}{{\bf D}}
\newcommand{\Xv}{{\bf X}}
\newcommand{\Yv}{{\bf Y}}
\newcommand{\Zv}{{\bf Z}}
\newcommand{\Uv}{{\bf U}}

\newcommand{\Vv}{{\bf V}}
\newcommand{\Rv}{{\bf R}}
\newcommand{\Qv}{{\bf Q}}

\newcommand{\Av}{{\bf A}}

\newcommand{\Pv}{{\bf P}}

\newcommand{\Sv}{{\bf S}}
\newcommand{\Nv}{{\bf N}}
\newcommand{\Iv}{{\bf I}}

\newcommand{\xv}{{\bf x}}

\newcommand{\tv}{{\bf t}}

\newcommand{\yv}{{\bf y}}

\newcommand{\zv}{{\bf z}}
\newcommand{\uv}{{\bf u}}

\newcommand{\nv}{{\bf n}}

\newcommand{\muv}{\boldsymbol \mu}

\newcommand{\Lambdav}{\boldsymbol \Lambda}

\newcommand{\Piv}{\boldsymbol \Pi}

\DeclareMathOperator\E{E}

 \def\E{\mathbb{E}}
 
 \def\Pr{\mathrm{Pr}}

\def\de \mathrm{d}

\newcommand{\Norm}{\mathcal{N}}

\newcommand\ie{i.e.,\xspace}
\def\textiid{i.i.d.\@\xspace}

\newcommand\iid{\ifmmode\text{ i.i.d. } \else \textiid \fi}

\newcommand{\Real}{\mathbb{R}}

\newcommand{\beqs}{\begin{equation*}}
\newcommand{\eeqs}{\end{equation*}}
\newcommand{\beq}{\begin{equation}}
\newcommand{\eeq}{\end{equation}}
\begin{document}
	\providecommand{\keywords}[1]{\textbf{\textit{Index terms---}} #1}

	\title{Shuffled Linear Regression via Spectral Matching}

	\IEEEoverridecommandlockouts
	\author{
		Hang Liu,~\IEEEmembership{Member,~IEEE}, and Anna~Scaglione,~\IEEEmembership{Fellow,~IEEE}
		\thanks{
		Hang Liu and Anna Scaglione are with the Department of Electrical and Computer Engineering, Cornell Tech, Cornell University, New York, NY, 10044 USA (e-mails: \{hl2382, as337\}@cornell.edu).
		}
	}
	\maketitle
\begin{abstract}
Shuffled linear regression (SLR) seeks to estimate latent features through a linear transformation, complicated by unknown permutations in the measurement dimensions. This problem extends traditional least-squares (LS) and Least Absolute Shrinkage and Selection Operator (LASSO) approaches by jointly estimating the permutation, resulting in shuffled LS and shuffled LASSO formulations. Existing methods, constrained by the combinatorial complexity of permutation recovery, often address small-scale cases with limited measurements. In contrast, we focus on large-scale SLR, particularly suited for environments with abundant measurement samples. We propose a spectral matching method that efficiently resolves permutations by aligning spectral components of the measurement and feature covariances. Rigorous theoretical analyses demonstrate that our method achieves accurate estimates in both shuffled LS and shuffled LASSO settings, given a sufficient number of samples. Furthermore, we extend our approach to address simultaneous pose and correspondence estimation in image registration tasks. Experiments on synthetic datasets and real-world image registration scenarios show that our method outperforms existing algorithms in both estimation accuracy and registration performance.
\end{abstract}

\section{Introduction}
Classical linear regression aims to recover multivariate signals from a noisy linear transformation, traditionally assuming that the correspondence between hidden signals and observations is perfectly known \emph{a priori}. However, recent research has indicated that in many practical scenarios the observations are subject to an unknown permutation, significantly complicating the traditional linear regression formulation \cite{SLR_ref3}. 
The challenge of permutation recovery has spurred the development of several solutions for the so-called \emph{shuffled linear regression (SLR)} \cite{SLR_ref3}, also referred to as unlabeled sensing \cite{8301561}. SLR aims to simultaneously recover hidden features and the associated permutations of measurements from the following model:
\begin{align}\label{eq01}
    \yv_i &= \Piv^\star \Av\xv_i^\star+\nv_i, i=1,\cdots,t,
\end{align}
where $\yv_i\in\Real^m$ and $\xv_i^\star\in\Real^n$ are the $i$-th measurement vector and latent feature vector, respectively, $\Av\in\Real^{m\times n}$ is a known sensing matrix, $\Piv^\star\in\mathcal{P}_m$ is an unknown permutation matrix from the set of $m\times m$ permutation matrices $\mathcal{P}_m$, and $\nv_i$ is an additive measurement noise vector independent of $\xv_i^\star$. 

The model in \eqref{eq01} finds widespread use across various real-world applications. In data de-anonymization, attackers attempt to uncover hidden features observed through a linear transformation, intertwined with unknown measurement shuffling that a publisher performed before making the data public \cite{SLR_ref5}. In image processing, SLR parallels the challenges encountered in classical \emph{image registration} problems \cite{10.1145/146370.146374}, which entail pose and correspondence estimation \cite{softposit}. This process involves aligning multiple point sets that capture the same objects at different snapshots. The aim of image registration is to identify the rigid transformation (referred to as the pose) and the permutation (referred to as the correspondence) between these datasets \cite{Revisited}.
Another significant application domain is the linear inversion system, particularly under conditions of signal \emph{asynchronization}. For instance, signal sampling with jitter \cite{1057717} illustrates a scenario where unknown observation times lead to unpredictably permuted measurements. This concern is particularly relevant in multi-target tracking and localization, which frequently involves estimating user positions from signals collected by asynchronous sensors \cite{POORE20061074}. Similar challenges arise in molecular channels \cite{7248460}, where receptors receive molecular tokens at varied intervals, each bearing indistinguishable signatures. 

Mathematically, for hidden features spanning real spaces as outlined in \eqref{eq01}, the tasks of simultaneous signal recovery and permutation estimation can be approached through a variant of the \emph{least-squares (LS)} method, often termed shuffled LS \cite{SLR_ref3}. Furthermore, exploiting the structural properties of the hidden features can enhance the estimation accuracy. Notably, the sparsity of the input signals proves beneficial in solving under-determined linear inversion problems within the compressed sensing framework. A prominent method in SLR involves adapting \emph{LASSO} algorithms to recover sparse signals while accommodating permutation estimation \cite{SLR_CS3}. In the subsequent section, we will review related work addressing these aspects.

\subsection{Related Work}
The SLR problem was initially explored in \cite{7447086}, where the authors established a necessary condition for exact permutation recovery in noiseless settings. The statistical limit for precise permutation recovery using an ideal maximum-likelihood estimator in a single-measurement scenario was subsequently investigated in \cite{SLR_ref3}. In a related context, Ref. \cite{SLR_ref7} analyzed the exact recovery conditions for shuffled LS algorithms with Gaussian sensing matrices. A branch-and-bound algorithm to solve shuffled LS was introduced in \cite{pmlr-v97-tsakiris19a} and later enhanced by concave minimization techniques in \cite{SLR_ref6}. Additionally, first-order methods aimed at enhancing computational efficiency were discussed in \cite{SLR_ref4b,SLR_ref8}. Research in more general multiple-measurement scenarios, corresponding to $t>1$ in \eqref{eq01}, was pursued in \cite{SLR_ref2} for maximum-likelihood estimation. Furthermore, Ref. \cite{pmlr-v119-zhang20n} built theoretical bounds for permutation estimation error under Gaussian sensing matrices and a restricted number of permutations.

In the context of sparse signal recovery, the authors in \cite{SLR_CS3} developed a branch-and-bound algorithm for solving shuffled LASSO. To reduce the computational complexity, a gradient-based approach for estimating multiple-measurement vectors and unknown permutations was introduced in \cite{SLR_CS2}. Additionally, Ref. \cite{ Alternating_Mat_SLR} proposed a permutation estimation algorithm based on spatial correlation matching and examined the impact of sparsity on estimation accuracy.

In addition to the standard SLR frameworks, various specialized variants have been developed to tackle different specific scenarios. A special case known as unlabeled de-noising, where $\Av$ is set to an identity matrix in \eqref{eq01}, was explored in \cite{SLR_ref5}. An alternative strategy focusing on unlabeled sub-sampling was analyzed in \cite{8234626}. Moreover, Ref. \cite{SLR_ref1} introduced a shuffled total LS method tailored for situations with noisy observations of the sensing matrix. Additionally, research in \cite{10.1214/18-EJS1498, SLR_SPermu_1} studied SLR problems involving sparse permutations, particularly under the assumption that only a small number of unknown permutations are possible in $\Piv^\star$.

Image registration, a common application of SLR, is a fundamental problem in image processing. The challenges of estimating pose (rigid transformation) and correspondence (permutation) have been extensively studied, but typically in isolation, each with a long-established history. Notably, LS methods for pose estimation were discussed in-depth in \cite{tpami,eggert1997estimating}, while spectral matching techniques for correspondence estimation were explored in \cite{1544893}. The combined challenge of simultaneous pose and correspondence estimation was first tackled in \cite{softposit}. Since this initial exploration, various algorithms have been developed to address the joint estimation problems, including branch-and-bound methods \cite{Revisited}, graduated non-convexity methods \cite{FGR}, and semi-definite relaxations \cite{9286491}, among others.

\subsection{Our Contributions}
The aforementioned work primarily focuses on scenarios with single measurements, corresponding to $t=1$ in \eqref{eq01}, or few-shot measurements, where $t$ is comparable to other parameters such as $m$. Typically, these methods search for the permutation that aligns observations and hidden features within the spatial domain. In contrast, the isolated permutation recovery challenges have been extensively addressed through spectral matching approaches that aim to align the eigenbases of the input and output data \cite{GMP_eign1}. These methods are acknowledged for their computational efficiency and robustness against noise, particularly when a substantial number of measurement samples are available \cite{GMP_eign3}.

Motivated by these considerations, in this work, we explore spectral matching to solve the SLR problems under the large-system limit. We consider the setup where the number of measurements $t$ significantly exceeds the dimension of measurements $m$, { i.e., $t\gg m$}. Our approach utilizes measurement samples to estimate the eigenbases of the signal sample covariance matrix. We then determine the optimal permutation by aligning the spectral components through spectral matching. This alignment allows us to reformulate the problem into traditional least squares (LS) or LASSO optimization frameworks. Finally, we extend our method to address the practical challenge of three-dimensional (3D) image registration. The key contributions of our work are summarized as follows:
\begin{itemize}
    \item We propose spectral matching methods for shuffled LS and shuffled LASSO problems and 
    { demonstrate} a reduced computational complexity compared to existing algorithms.
    \item We show that the permutation estimation error of our proposed method diminishes at the rate of $\mathcal{O}(1/\sqrt{t})$ when $t \gg m^3$. By comparing to the state-of-the-art algorithm \cite{pmlr-v119-zhang20n}, which computes the permutation in the spatial domain, our analysis underscores the superior performance of spectral matching, particularly when a large number of observation samples are available.
    \item Based on permutation error analysis, we study the estimation error for both shuffled LS and LASSO. Our method achieves sub-linear error rates, balancing efficiency with optimal asymptotic order as $t \to \infty$.  Specifically, as $t \to \infty$, the estimation error of our method achieves the optimal asymptotic order, albeit with a slower convergence rate ($\mathcal{O}(1/\sqrt{t})$) compared to the optimal but computationally prohibitive solver whose rate is $\mathcal{O}(1/t)$. 
    \item We adapt the proposed method for 3D image registration by combining spectral matching with a translation estimation step. While image registration does not typically follow the large-system limit, the orthogonality of the rotation matrix enables accurate covariance estimation, broadening our method’s applicability.
\end{itemize}

We conduct experiments on synthetic data and real-world 3D point datasets for addressing both shuffled LS and shuffled LASSO problems. The results demonstrate that our proposed method surpasses state-of-the-art baselines in performance and achieves accurate permutation estimation, particularly when the sample size $t$ is sufficiently large.

\subsection{Organization and Notations}
The remainder of this paper is organized as follows. We formulate the SLR problems in Section \ref{sec2}. In Section \ref{sec3}, we develop the spectral matching algorithm.  In Section \ref{sec4}, we analyze the estimation error of the proposed method. In Section \ref{sec5}, we adapt the proposed method to solve 3D image registration problems. In Section \ref{sec6}, we present numerical results to evaluate the proposed method. Finally, this paper concludes in Section \ref{sec7}.
	
	Throughout, 
	we use regular, bold small, and bold capital letters to denote scalars, vectors, and matrices, respectively. We use  $\Xv^T$ to denote the transpose of matrix $\Xv$, $\tr(\Xv)$ to denote its trace, and $\rank(\Xv)$ to denote its rank.
	We use $x_i$ to denote the $i$-th entry of vector $\xv$, $x_{ij}$ or $[\Xv]_{i,j}$ interchangeably to denote the $(i,j)$-th entry of matrix $\Xv$, and $\xv_j$ to denote the $j$-th column of $\Xv$. The real normal distribution with mean $\muv$ and covariance $\Cv$ is denoted by $\Norm(\muv,\Cv)$, and the cardinality of set $\mathcal{S}$ is denoted by $\abs{\mathcal{S}}$. We use $\norm{\cdot}_p$ to denote the $\ell_p$ norm,  $\norm{\cdot}_F$ (resp. $\norm{\cdot}_2$) to denote the Frobenius (resp. spectral) matrix norm, $\Iv_n$ to denote the $n\times n$ identity matrix,  $\bf 1$ (resp. $\bf 0$) to denote the all-one (resp. all-zero) vector with an appropriate size, and $\diag(\xv)$ to denote a diagonal matrix with the diagonal entries specified by $\xv$. For two scalars $a$ and $b$, we denote $a\gtrsim b$ or $b\lesssim a$ interchangeably if there exists a positive constant $c_0$ such that $a\geq c_0b$. We denote $a \asymp b$ when both $a\gtrsim b$ and $b\gtrsim a$ hold.

\section{Problem Formulation}\label{sec2}
We study the SLR model described in \eqref{eq01}. Our goal is to estimate both the permutation $\Piv^\star$ and the latent features $\Xv^\star=[\xv_1^\star,\cdots,\xv_t^\star]$ from the noisy measurements $\Yv=[\yv_1,\yv_2,\cdots,\yv_t]$ and the sensing matrix $\Av$. Unless otherwise specified, the feature vectors $\{\xv_i^\star\}_{i=1}^t$ and the noise vectors $\{\nv_i\}_{i=1}^t$ are modeled as unknown independent and identically distributed (i.i.d.) random vectors drawn from unspecified distributions. These vectors are assumed to have zero means with given positive semidefinite covariance matrices $\E[\xv_i^\star(\xv_i^\star)^T]=\Cv_X$ and $\E[\nv_i\nv_i^T]=\Cv_N$. The adaptation of the proposed estimation approach to accommodate features with non-zero means and certain dependent samples will be explored in Section \ref{sec5}.

Under the condition of $m\geq n$, the SLR estimation problem has been explored through a \emph{shuffled LS} formulation \cite{SLR_ref6}:
\begin{align}\label{eq02}
    (\widehat \Piv,\widehat\Xv)\in\argmin_{\Xv\in\Real^{n\times t},\Piv\in\mathcal{P}_m}\norm{\Piv^T\Yv-\Av\Xv}_F^2.
\end{align}
It has been shown that the problem in \eqref{eq02} is NP-hard \cite{SLR_ref3}.

In addition to the shuffled LS, we examine a variant of this problem that accommodates underdetermined linear regression with $n \geq m$, specifically when the input features in $\Xv$ are \emph{sparse}. This setup generalizes the well-established compressed sensing framework \cite{1614066} with the additional complexity of an unknown measurement permutation. The sparsity assumption models scenarios where signals are sparsely represented under certain linear transformations, such as image pixels under wavelet or Fourier transforms. Previous work \cite{SLR_CS3, Alternating_Mat_SLR} proposed to solve compressed sensing under unknown permutations using the following $\ell_1$-regularized formulation:
\begin{align}\label{eq03}
    (\widehat \Piv, \widehat \Xv) \in \argmin_{\Xv \in \Real^{d \times t}, \Piv \in \mathcal{P}_m} \norm{\Piv^T \Yv - \Av \Xv}_F^2 + \rho \norm{\Xv}_{1,1},
\end{align}
where $\rho$ is a tunable penalty coefficient and $\norm{\Xv}_{1,1}=\sum_{i,j}x_{ij}$ denotes the entrywise $\ell_1$-norm. We refer to \eqref{eq03} as \emph{shuffled LASSO} regression.

The goal of this work is to design computationally efficient algorithms for solving \eqref{eq02} and \eqref{eq03}, for cases where a substantial number of measurement samples is available, i.e., $t \gg m$. 

\remark{We note that the problem in \eqref{eq03} is fundamentally different from those considered in \cite{10.1214/18-EJS1498,SLR_SPermu_1}. The latter studies typically involve a dense feature matrix $\Xv$ and a sparse permutation $\Piv^\star$. In contrast, our formulation does not impose any restrictions on the unknown permutation. Consequently, in the worst-case scenario, there could be $m!$ possible random permutations to search.}

\section{Spectral Matching}\label{sec3}
In this section, we present an approach, with relatively low-complexity, to solve \eqref{eq02} and \eqref{eq03}, which is particularly suitable for the high-dimensional scenarios with $t \gg m$. The main challenge in these combinatorial problems stems from the unknown permutation matrix, casting the issue as a quadratic assignment problem, as discussed in \cite{SLR_CS2}. Once a solution for the permutation $\widehat\Piv$ is identified, the solution to $\Xv$ in \eqref{eq02} exhibits a closed-form expression as
\begin{align}\label{eq04}
\widehat \Xv(\widehat\Piv) = \Av^\dagger \widehat\Piv^T \Yv,
\end{align}
where $\Av^\dagger = (\Av^T\Av)^{-1}\Av^T$ is the pseudo-inverse of $\Av$. With the LASSO formulation, once $\widehat\Piv$ is determined, problem \eqref{eq03} can be mapped into:
\begin{align}\label{eq04b}
\widehat \Xv(\widehat\Piv) = \argmin_{\Xv \in \Real^{d \times t}} \norm{\widehat\Piv^T \Yv - \Av \Xv}F^2 + \rho \norm{\Xv}_{1,1},
\end{align}
which is a convex programming problem and can be efficiently solved using off-the-shelf optimization solvers.

In the subsequent analysis, we focus on solving for $\widehat\Piv$ from the model in \eqref{eq01}. Our approach hinges on estimating the covariance matrix of $\Yv$ and aligning its eigenvectors with those of $\Cv_X$. Specifically, we denote the covariance of $\yv_i$ by $\Cv_Y = \E[\yv_i \yv_i^T]$. From \eqref{eq01}, the covariance is expressed as:
\begin{align}
    \Cv_Y = \Piv^\star \Av \Cv_X \Av^T (\Piv^\star)^T + \Cv_N.
\end{align}
Consequently, $\Cv_Y -\Cv_N$ equals to $\Av \Cv_X \Av^T$ subject to the permutation of $\Piv^\star$ applied to its rows and columns. We detail their eigendecompositions as:
\begin{align}
    \Cv_Y - \Cv_N &= (\Piv^\star \Vv) \Lambdav (\Piv^\star \Vv)^T, \nonumber\\
    \Av \Cv_X \Av^T &= \Vv \Lambdav \Vv^T,\label{eq06}
\end{align}
where $\Lambdav = \diag(\lambda_1, \cdots, \lambda_m)$ is a diagonal matrix containing the eigenvalues in descending order, and $\Vv$ is an orthogonal eigenmatrix of $\Av \Cv_X \Av^T$. 

As a result, the eigenvectors of $\Cv_Y - \Cv_N$ correspond to those of $\Av \Cv_X \Av^T$, subject to the permutation $\Piv^\star$. Therefore, $\Piv^\star$ can be determined by aligning their eigenvector matrices. To facilitate this, we compute the sample covariance matrix of $\Yv$, denoted by $\widehat \Cv_Y$, to approximate $\Cv_Y$ and represent the corresponding eigendecomposition as:
\begin{align}\label{eq07}
    \widehat \Cv_Y - \Cv_N = \frac{1}{t} \sum_{i=1}^t \yv_i \yv_i^T - \Cv_N = \Uv \Sv \widehat \Lambdav \Sv \Uv^T,
\end{align}
where $\widehat \Lambdav = \diag(\hat \lambda_1, \cdots, \hat \lambda_m)$ consists of the sample eigenvalues in descending order, $\Uv$ is the sample eigenvector matrix, and $\Sv$ is a diagonal matrix with entries of either $1$ or $-1$. The matrix $\Sv$ accounts for the inherent sign ambiguity in eigendecompositions.\footnote{If $(\uv,\lambda)$ is an eigenpair for a matrix, then $(-\uv, \lambda)$ is also a valid eigenpair for the same matrix.}

With a large $t$, the sample covariance $\widehat \Cv_Y$ provides a good approximation of $\Cv_Y$. By combining \eqref{eq06} and \eqref{eq07}, we see that $\Uv$ closely approximates $\Piv^\star\Vv\Sv$. This motivates us to compute the permutation by minimizing the distance between  $\Uv$ and $\Piv^\star\Vv\Sv$. To tackle the sign ambiguity, we compute the permutation by solving the following problem (cf. \cite{GMP_eign1,BlindGM}):
\begin{align}\label{eq08}
    \widehat\Piv=\argmax_{\Piv\in\mathcal{P}_m}\tr\left( abs\left(\Vv_{\mathcal{A}}\right)abs\left(\Uv^T_{\mathcal{A}}\right)\Piv \right),
\end{align}
where ${abs}(\cdot)$ is the matrix containing the absolute values of the entries in the input matrix, and $\Vv_{\mathcal{A}}$ denotes the submatrix consisting of the columns of $\Vv$ indexed by a subset $\mathcal{A} \subset \{1, 2, \cdots, m\}$. Problem \eqref{eq08} is a linear assignment problem and can be efficiently solved using the Hungarian method \cite{Hungarian}.

To ensure precise alignment of eigenvectors, we propose to select $\mathcal{A}$ from the \emph{principal components} of $\Vv$ that correspond to distinct and non-zero eigenvalues as
\begin{align}\label{eq10}
    \mathcal{A}\subset \mathcal{C}\triangleq\{1\leq i\leq m:\lambda_i> 0,\lambda_i\neq \lambda_j, \forall j\neq i\}.
\end{align}
Given that $\rank(\Av\Cv_X\Av^T)\leq \min\{m,n\}$, we have $|\mathcal{A}|\leq |\mathcal{C}|\leq \min\{m,n\}$. Although one could set $|\mathcal{A}| = |\mathcal{C}|$, it has been shown in \cite{BlindGM} that selecting only a subset of principal components, particularly those corresponding to large spectral gaps, boosts estimation accuracy. For instance, $\mathcal{A}$ can be progressively computed using the line-search algorithm in \cite[Algorithm 3]{BlindGM}. The selection of $\mathcal{A}$ from the principal components guarantees a strictly positive spectral gap among the chosen spectral components. Specifically, for any given $\mathcal{A}$, the minimum spectral gap is denoted by: 
\begin{align} 
\delta\triangleq \min_{i\in\mathcal{A}} \{\lambda_i-\lambda_{i+1}\}, 
\end{align} 
where $\lambda_0=\infty$ and $\lambda_{m+1}=0$ are defined for notational convenience. The uniqueness of the selected eigenvalues ensures that $\delta > 0$. This condition is crucial for ensuring the correct alignment of the sample eigenvectors, which is a prerequisite for achieving accurate permutation estimation.

    \begin{figure}[!t]
                 	\centering
                 	\includegraphics[width=3.2 in]{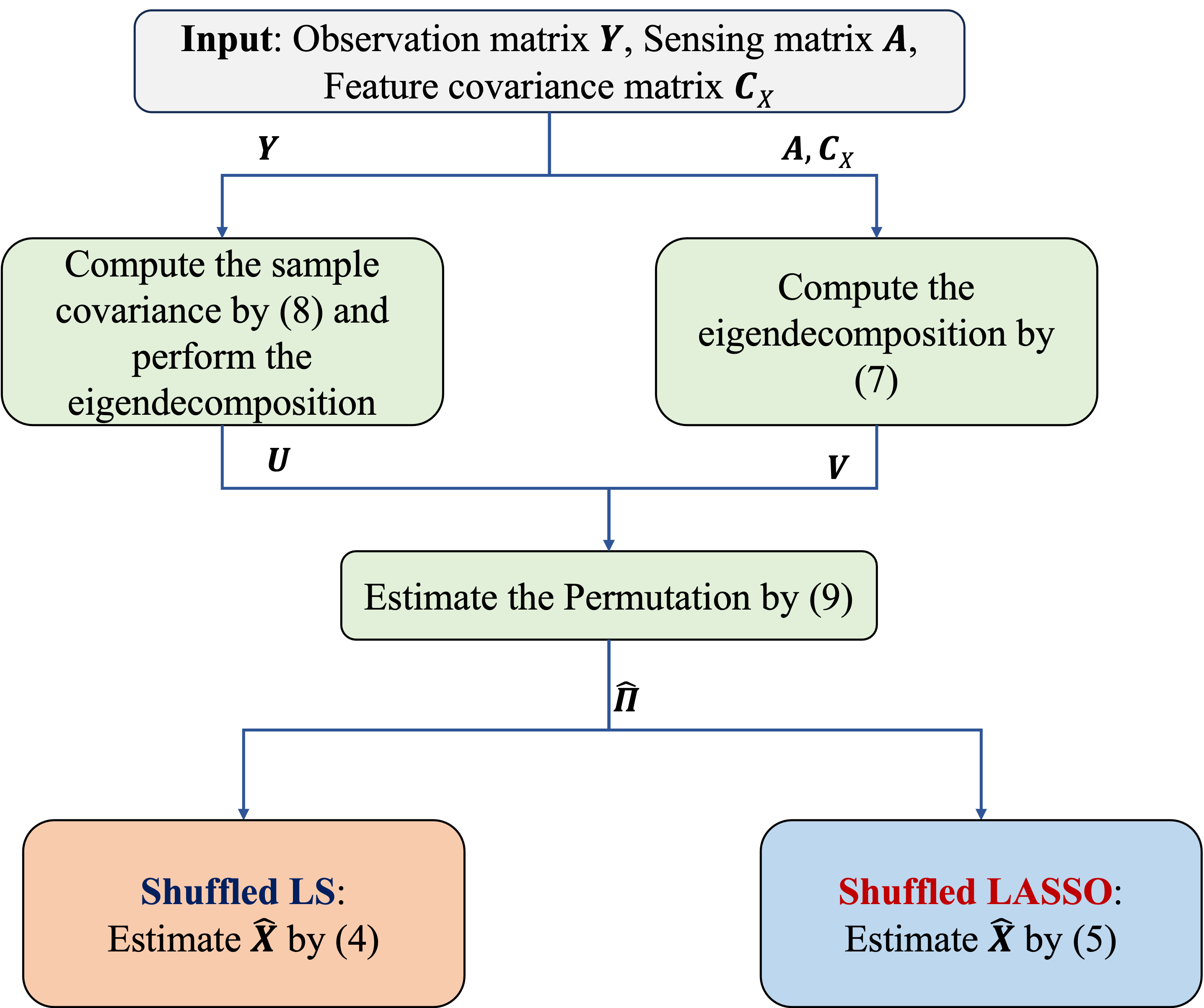}
                 	\caption{Schematic view of the spectral matching method.}
                 	\label{diagram}
    \end{figure}
                 
In summary, our approach first computes the sample covariance matrix and its eigendecomposition by \eqref{eq06} and \eqref{eq07}. Subsequently, we estimate $\widehat\Piv$ by solving the linear assignment problem in \eqref{eq08}. The solution for $\Xv$ is then derived either directly from \eqref{eq04} or by solving the resultant convex problem in \eqref{eq04b}. { The whole spectral matching process is summarized in Fig. \ref{diagram}.}

The total computational complexity for computing \eqref{eq06}--\eqref{eq08} scales at the rate of $\mathcal{O}(m^2(m+t))$. In contrast, the state-of-the-art approximate branch-and-bound algorithm in \cite{SLR_ref6}, applicable only when $t=1$, exhibits a prohibitive complexity of $\mathcal{O}(n^{2m})$ for large values of $n$. Additionally, the complexity of the spatial-domain matching algorithm in \cite{pmlr-v119-zhang20n} grows at a rate of $\mathcal{O}(m^2(m+n+t))$.

{ \remark{We note that a strictly positive spectral gap $\delta$ in the selected eigencomponents is crucial for the effectiveness of our proposed approach. To achieve this, we prioritize the selection of principal eigencomponents that correspond to distinct eigenvalues in \eqref{eq10}. This means that, in principle, our method works only if the matrix $\Av \Cv_X \Av$ has at least one eigencomponent corresponding to distinct and positive eigenvalues. In other words, the number of selected eigenvectors should be at least one, i.e., $\mathcal{A} \neq \emptyset$. We do acknowledge an extreme case that would violate this assumption: if all the eigenvalues of $\Av \Cv_X \Av$ are identical, i.e., when $\Av \Cv_X \Av$ is proportional to an identity matrix $\Iv_m$. In this case, it is impossible to distinguish the order from the spectrum.}}
\section{Performance Analysis}\label{sec4}
We analyze the estimation error for the proposed solutions to shuffled LS and shuffled LASSO. As detailed in Section \ref{sec2}, we focus on the regime where $t \gg m$.

\subsection{Permutation Estimation Error Analysis}
We begin by examining the permutation estimation error for the solution in \eqref{eq08}. { We denote the number of selected eigenvectors in the considered linear assignment problem in \eqref{eq08} as 
\begin{align}
    k\triangleq |\mathcal{A}|\leq \min\{m,n\}.
\end{align} 
}

To assess the accuracy of the estimate $\widehat\Piv$, we introduce a noiseless permutation estimation problem as a counterpart to \eqref{eq08}:
\begin{align}\label{eq_noiselesspro}
    \max_{\Piv\in\mathcal{P}_m}\tr\left(\Piv^T\Piv^\star \Vv_{\mathcal{A}}\Vv^T_{\mathcal{A}} \right),
\end{align}
where the sign ambiguity no longer exists, thus negating the need for the absolute value operation. Problem \eqref{eq_noiselesspro} mirrors the ideal scenario where the covariance estimation in \eqref{eq07} is perfectly noiseless, \ie $\widehat\Cv_Y=\Cv_Y$. We note that \eqref{eq_noiselesspro} is an oracle formulation as it presupposes the knowledge of $\Piv^\star$. It can be verified that for any $\Piv\in\mathcal{P}_m$,
\begin{align}\label{eq012}
    \tr\left(\Piv^T\Piv^\star \Vv_{\mathcal{A}}\Vv^T_{\mathcal{A}} \right)\leq \tr\left( \Vv_{\mathcal{A}}\Vv^T_{\mathcal{A}}\right)=k,
\end{align}
where equality holds if $\Piv=\Piv^\star$. In other words, the true permutation $\Piv^\star$ maximizes the objective in \eqref{eq_noiselesspro}. 
Hence a measure of the permutation accuracy of our solution $\widehat\Piv$ is the normalized \emph{optimality gap} with respect to \eqref{eq_noiselesspro}, defined as
\begin{align}\label{eq13}
    1-\frac{1}{k}\tr\left(\widehat\Piv^T\Piv^\star \Vv_{\mathcal{A}}\Vv^T_{\mathcal{A}} \right)\in[0,1].
\end{align}
 Intuitively, \eqref{eq13} quantifies the average mismatched ratio for the proposed solution, which equals to zero if $\widehat\Piv=\Piv^\star$. The following theorem bounds this optimality gap. 
\theorem[Optimality gap]{\label{theorem1} { Suppose that the number of selected eigenvectors in \eqref{eq10} is at least one, i.e., $k\geq 1$.}
Moreover, suppose the $\ell_2$-norm of each measurement vector $\norm{\yv_i}_2$ in \eqref{eq01} is uniformly bounded from above almost surely. Let $\mathcal{A}$ be defined according to \eqref{eq10} with $k=|\mathcal{A}|\leq\min\{m,n\}$. For any $\epsilon>0$ and $t\geq \frac{m(\ln m+\epsilon)}{\min\{1,\delta^2\}}$, with probability at least $1-2e^{-\epsilon}$, the solution in \eqref{eq08} satisfies
\begin{align}\label{eq14}
 1- \frac{1}{k}\tr(\widehat\Piv^T\Piv^\star \Vv_{\mathcal{A}}\Vv^T_{\mathcal{A}})\lesssim \frac{m}{\delta}\sqrt{\frac{m(\ln m +\epsilon)}{kt}}.
\end{align}
}
\begin{IEEEproof}
See Appendix \ref{appa}.
\end{IEEEproof}

For a fixed $m$, Theorem \ref{theorem1} shows that the average mismatch ratio decays at the rate of $\mathcal{O}(1/\sqrt{t})$ provided that $t\gtrsim m\ln m$. We note that Theorem \ref{theorem1} does not require stringent assumptions on the distributions of the measurement and noise vectors.

In addition to assessing the optimality gap, we investigate the probability that the proposed solution incurs at least one mismatch relative to the true permutation. To this end, we impose an additional assumption regarding the distribution of the hidden features as follows.

\definition[Rotational invariant distributions]{The probability distribution of a random vector $\xv\in\Real^m$ is considered rotationally invariant if, for any rotation matrix $\Rv\in SO(m)=\{\Rv^T\Rv=\Iv_m,|\det(\Rv)|=1\}$, we have $\Pr(\xv)=\Pr(\Rv\xv)$.
}

Rotational invariance implies that the distribution is unchanged under rotations. This property is satisfied with a broad range of widely used distributions in physics, signal processing, and machine learning. Common examples include Gaussian distributions, Student's t-distributions, and uniform distributions on the sphere. The following theorem explores the error probability associated with rotationally invariant feature vectors.

\theorem[Permutation error probability]{\label{theorem2}
Suppose the following conditions hold:
\begin{enumerate}
    \item[i.] Each $\{\yv_i\}_{i=1}^m$ is uniformly bounded above almost surely.
    \item[ii.] The distribution of $\xv_i^\star$ is rotational invariant. 
    \item[iii.] $t\gtrsim \frac{mk\ln m}{\delta^2}$.
    \item[iv.] { $k\geq 1$.} 
\end{enumerate}
Then, we have
\begin{align}\label{eq15}
    \Pr(\widehat\Piv\neq \Piv^\star)\lesssim\sqrt{\frac{m^3k^5}{m^3/\ln m+\delta^2t}}e^{-c_0\frac{kt}{\ln m(m^3/\delta^2+t)}},
\end{align}
where $c_0$ is an absolute constant independent of $m,n,k,$ and $t$.
Furthermore, if $k\gtrsim \ln m$ and $t\gtrsim m^3 k^5/\delta^2$, there exists a constant $t_0$ such that for any $t\geq t_0$ we have
\begin{align}\label{eq16}
    \Pr(\widehat\Piv\neq \Piv^\star)\lesssim \frac{mk^2}{\delta}\sqrt{\frac{mk}{t}}\overset{t\to\infty}{\longrightarrow}0.
\end{align}
}
\begin{IEEEproof}
See Appendix \ref{appb}.
\end{IEEEproof}

Theorem \ref{theorem2} shows that the error probability diminishes at the rate of $\mathcal{O}(1/\sqrt{t})$ provided that $t\gtrsim m^3(\ln m)^5$. We emphasize that this diminishing rate holds even with finite values of $m$ and $k$ as long as the measurement sample size $t$ is sufficiently large. Such a diminishing error rate underscores the advantages of \emph{spectral}-domain matching in \eqref{eq08}. Specifically, with sufficiently many measurement samples, matching permutations based on sample eigenvectors proves to be more robust against measurement noise than existing spatial-domain matching methods. To verify the claimed benefits of spectral matching, we compare our results with a state-of-the-art spatial-domain matching approach in \cite[Eq. (4)]{pmlr-v119-zhang20n}:
    \begin{align}\label{eq17}
\widehat\Piv_{\text{sp}}=\argmax_{\Piv\in\mathcal{P}_m}\tr\left(\Piv\Av\Av^T\Yv\Yv^T\right).
    \end{align}
Intuitively, this approach relies on the expectation that $\Av^T\Yv$ approximates $\Xv^\star$, given near-orthogonal projection and minimal noise. Therefore, the formulation in \eqref{eq17} essentially seeks to maximize the spatial correlation of the measurements and the features. The following lemma states the achievable error probability for \eqref{eq17}. 
\lemma[Achievable error bound in \cite{pmlr-v119-zhang20n}]{\label{lemma1}Suppose the following conditions hold:
\begin{enumerate}
    \item[i.] The entries of $\Av$ are i.i.d. drawn from a standard Gaussian distribution $\Norm(0,1)$. Moreover, each noise vector $\nv_i$ follows the distribution of $\Norm({\bf 0},\sigma^2\Iv_m)$.
    \item[ii.] $m\gtrsim n^4(\log n)^4$.
    \item[iii.] The number of non-zero off-diagonal entries in $\Piv^\star$ is less than $0.25m$, i.e.,  $\tr(\Piv^\star)\geq 0.75m$.
\end{enumerate}
Under these conditions and assuming that $m,n,$ and $\norm{\Xv^\star}_F^2$ are sufficiently large, the solution in \eqref{eq17} satisfies
 \begin{align}\label{lemma1result}
     \Pr(\widehat\Piv_{\text{sp}}\neq \Piv^\star)\lesssim me^{-c_o^\prime\cdot\min\{m,t\}}+\frac{c_1^\prime}{n^2},
 \end{align}
where $c_0^\prime$ and $c_1^\prime$ are some constants.
}
\begin{IEEEproof}
    See Appendix \ref{app_lemma1}.
\end{IEEEproof}

Comparing Theorem \ref{theorem2} with Lemma \ref{lemma1}, our results are achieved without strict assumptions on the distributions of $\Av$ and $\nv_i$, nor on the sparsity of the true permutations. Furthermore, when $m$ and $n$ are fixed and $t \geq m$, the error probability for \eqref{eq17} is bounded by $\mathcal{O}(me^{-m}+1/n^2)$, which does not improve as $t$ increases. This limitation stems from the fact that $\frac{1}{t} \Yv \Yv^T$ stabilizes for large $t$; hence the error in \eqref{eq17} is primarily due to misalignment in $\Av^T \Yv$. Conversely, our spectral matching approach benefits from a larger sample size, enhancing eigenvector estimation accuracy. This leads to a progressively diminishing error probability as $t$ increases.

{ \remark{Theorems \ref{theorem1} and \ref{theorem2} emphasize the significant role of $\delta$ in determining the accuracy of the matching process. Specifically, the permutation error is inversely proportional to $\delta$. Intuitively, when the spectral gap is small, observation noise is more likely to disrupt the order of adjacent sample eigenvectors, leading to incorrect permutation alignment. This aligns with conventional perturbation theory, which indicates that a smaller spectral gap increases the sensitivity of spectrum estimation to noise perturbations. As we demonstrate in the subsequent section, this misalignment in the permutation estimation propagates, affecting the estimation accuracy of $\Xv$.
}
}

\subsection{Shuffled LS Performance}
Building on the permutation error bounds, we investigate the accuracy of the proposed solution in \eqref{eq04} and \eqref{eq08} applied to the shuffled LS problem in \eqref{eq02}. We examine a deterministic system in \eqref{eq01} with the condition $m \geq n$ and Gaussian measurement noise. To establish a baseline for performance evaluation, we first reference the optimal reconstruction mean-square error (MSE) achievable by the \emph{best possible} shuffled LS solution, as stated in \cite[Theorem 1]{SLR_ref2}:
\lemma[Optimal reconstruction MSE \cite{SLR_ref2}]{\label{lemma2}Suppose that the measurement noise in \eqref{eq01} follows the Gaussian distribution of $\Norm({\bf 0},\sigma^2\Iv_m)$. Let $m\geq n$ and $\rank(\Av)=n$. Denote the optimal solution to \eqref{eq02} by $(\widehat\Piv_{\text{opt}},\widehat\Xv_{\text{opt}})$. With probability at least $1-e^{-Const\cdot (m\ln m+nt)}$,
\begin{align}
   \frac{ \norm{\widehat\Piv_{\text{opt}}\Av\widehat\Xv_{\text{opt}}-\Piv^\star\Av\Xv^\star}_F^2}{mt}\lesssim \sigma^2\left(\frac{n}{m}+\frac{\ln m}{t}\right).
\end{align}
}

As $t\to \infty$, Lemma \ref{lemma2} shows that the optimal solution to the shuffled LS problem in \eqref{eq02} converges to the same reconstruction error rate as the solution to conventional LS problems, specifically $\mathcal{O}(\frac{n}{m}\sigma^2)$. However, the shuffled LS problem in \eqref{eq02} is  NP-hard \cite{SLR_ref2}, making the optimal solution generally intractable. We refer to the result in Lemma \ref{lemma2} as the optimal error rate for shuffled LS. Following this, we present the reconstruction error bound of our proposed approach.
\theorem{\label{theorem3} Suppose 1) { $k\geq 1$}, 2) each $\yv_i$ in \eqref{eq01} has a uniformly bounded norm almost surely, and 3) $\nv_i$ follows the Gaussian distribution of $\Norm({\bf 0},\sigma^2\Iv_m)$. Denote the proposed solution to \eqref{eq02} in \eqref{eq04} and \eqref{eq08} by $\widehat\Xv$ and $\widehat\Piv$. When $m\geq n$, $t\gtrsim \frac{m\ln m}{\min\{1,\delta^2\}}$, and $k\gtrsim n$, with probability at least $1-\frac{4}{m}-e^{-nt/8}$, we have 
\begin{align}\label{eq20}
   \frac{ \norm{\widehat\Piv\Av\widehat\Xv-\Piv^\star\Av\Xv^\star}_F^2}{mt}\lesssim \frac{n}{m}\sigma^2+\frac{n}{\delta}\sqrt{\frac{\ln m}{t}}.
\end{align}
}
\begin{IEEEproof}
See Appendix \ref{appc}.
\end{IEEEproof}
As $t\to\infty$, the reconstruction error of the proposed method approaches $\mathcal{O}(\frac{n}{m}\sigma^2)$, converging at a sub-linear rate of $\mathcal{O}(1/\sqrt{t})$. Compared with Lemma \ref{lemma2}, our method reaches the asymptotically optimal error bound, albeit at a slower rate of convergence ($\mathcal{O}(1/\sqrt{t})$ compared to $\mathcal{O}(1/t)$). On the other hand, while the optimal estimator in Lemma \ref{lemma2} typically exhibits an exponential computational complexity, our method benefits from a polynomial complexity; see Section \ref{sec3}. Moreover, we note that the bound in \eqref{eq20} requires the number of selected eigenvectors $k$ to scale at least linearly with the dimension of the features $n$. This scaling is crucial for achieving accurate estimation as the system size increases.

\subsection{Shuffled LASSO Performance}
We turn to the performance analysis on the shuffled LASSO problem in \eqref{eq03}. Traditional compressed sensing theory highlights the utility of exploiting feature sparsity to resolve underdetermined problems even when $n > m$. With a sufficiently large number of measurements $t$, Theorem \ref{theorem2} illustrates that our permutation estimate converges to the true permutation, thus making classical compressed sensing results applicable. To illustrate this, we consider a scenario where each feature vector $\xv_i^\star$ is $s$-sparse, meaning the number of non-zero entries $\norm{\xv_i^\star}_0$ is at most $s$. Under this assumption, it is well-established that reliable recovery is feasible through the traditional LASSO estimator, provided that the sensing matrix $\Av$ satisfies the restricted eigenvalue condition; see \cite[Definition 7.12]{HDimPro2}. Integrating these classical LASSO results with the permutation error bound provided in Theorem \ref{theorem2}, we establish the following error bound for the shuffled LASSO problem.
\theorem{\label{theorem4} 
Consider the solution $\widehat\Xv$ and $\widehat\Piv$ in \eqref{eq04b} and \eqref{eq08} to the shuffled LASSO problem \eqref{eq03}. Suppose the following conditions hold:
\begin{enumerate}
    \item[i.] The conditions in Theorem \ref{theorem2} hold.
    \item[ii.] $k\gtrsim \ln m$.
    \item[iii.] The noise $\nv_i$ in \eqref{eq01} follows the distribution of $\Norm({\bf 0},\sigma^2\Iv_m)$.
    \item[iv.] The penalty coefficient in \eqref{eq03} satisfies $\rho\gtrsim \sigma \sqrt{\frac{\ln n}{m}}$.
    \item[v.] Each feature vector $\xv_i^\star$ is $s$-sparse, i.e., $\norm{\xv_i^\star}_0\leq s,\forall i$.
    \item[vi.] The sensing matrix $\Av$ satisfies the restricted eigenvalue condition over the support of every  $\xv_i^\star$ with parameter $(\kappa,3)$. Specifically, denoting the support of $\xv_i^\star$ by $\mathcal{S}_i$, it holds for any $1\leq i\leq t$ that \cite[Definition 7.12]{HDimPro2}
    \begin{align*}
        \frac{1}{m}\norm{\Av\zv}^2_2\geq \kappa \norm{\zv}_2^2, \forall \zv\in\left\{\zv:\sum_{k\notin\mathcal{S}_i}|z_j|\leq 3 \sum_{k\in\mathcal{S}_i}|z_j|\right\}.
    \end{align*}
\end{enumerate}

Then, with probability at least $1-\mathcal{O}(\frac{1}{n}+\sqrt{\frac{m^3k^5}{\delta^2t}})$, there exists a constant $t_0$ such that for any $t\geq t_0$ the following two bounds hold:
\begin{align}
      &\frac{ \norm{\widehat\Piv\Av\widehat\Xv-\Piv^\star\Av\Xv^\star}_F^2}{mt}\lesssim  \frac{s\ln n }{\kappa m}\sigma^2,\\
      &\frac{ \norm{\widehat\Xv-\Xv^\star}_F^2}{t}\lesssim  \frac{s\ln n }{\kappa m}\sigma^2.
\end{align}
}
\begin{IEEEproof}
   The result follows from the combination of Theorem \ref{theorem2} and \cite[Theorems 7.13 and 7.20]{HDimPro2}.
\end{IEEEproof}

Similar to Theorem \ref{theorem3}, the error bound for our approach converges to $\mathcal{O}(\frac{s\ln n}{m}\sigma^2)$ as $t\to\infty$. When $s\leq n/\ln n$, the shuffled LASSO estimator achieves a smaller error compared to shuffled LS due to the exploitation of sparsity.

\section{Applications to Image Registration}\label{sec5}
In this section, we study the application of the proposed method in image registration. Specifically, we aim to solve the simultaneous pose and correspondence estimation problem \cite{softposit}, which involves registering a 3D point set to a 2D image or another 3D point set without knowing point correspondence. Denote by $\Qv\in \Real^{m\times 3}$ the given 3D point set, where each row represents the Cartesian coordinates of each point. Consider the projection of a 3D point set onto a 2D image. We assume a transparent 3D model so that all the points preserve on the 2D projected image. The homogeneous coordinates corresponding to this projected image, which originate from the center of projection, are denoted by $\Pv \in \Real^{m \times 3}$. The relationship between $\Pv$ and $\Qv$ is given by \cite[Section 2]{softposit}
\begin{align}\label{eq_irmodel}
    \Pv=\Piv^\star\Qv\Rv+{\bf 1}\tv^T,
\end{align}
where $\Piv^\star\in\mathcal{P}_m$ is the unknown permutation matrix representing the point correspondence, $\Rv$ represents the unknown rotation matrix due to the coordinate transformation belonging to the rotation group $SO(3)\triangleq\{\Rv\in\Real^{3\times 3}:\Rv^T\Rv=\Iv_3,\det(\Rv)=1\}$, and $\tv\in\Real^3$ represents the unknown translation vector.

The parameters $\Piv^\star$, $\Rv$, and $\tv$ in \eqref{eq_irmodel} model the coordinate transformation of the projection, accompanied by the loss of point correspondence. The objective is to recover the pose ($\Rv$ and $\tv$) and the correspondence ($\Piv^\star$) from the (noisy) observations of $\Pv$ and $\Qv$. This challenge is thus termed the {\it simultaneous pose and correspondence problem}. For further details on this projection model, we refer to \cite[Sections 2-4]{softposit}. Additionally, the model in \eqref{eq_irmodel} also describes the 3D-3D registration problem, which involves determining the correspondence between two point sets under an unknown permutation and rigid transformation \cite{Revisited}.

Given the (potentially noisy) observations of $\Qv$ and $\Pv$, the problem can be solved via the shuffled LS model in \eqref{eq02} as  \cite{Revisited}
\begin{align}\label{eq_ir}
    (\widehat \Piv,\widehat\Rv,\hat \tv)\in\argmin_{\Rv\in SO(3),\Piv\in\mathcal{P}_m,\tv\in\Real^{3}}\norm{\Pv-\Piv\Qv\Rv-{\bf 1}\tv^T}_F^2.
\end{align}
The model in \eqref{eq_ir}, which can be analogously related to \eqref{eq02} with $t=3$, at first glance, does not appear to meet the condition of $t \gg m$ necessary for applying the spectral matching method described in Section \ref{sec3}. However, the intrinsic orthogonality of the rotation matrix $\Rv$ yields the following simplification:
\begin{align}
   \Pv\Pv^T=\Piv^\star\Qv\Qv^T(\Piv^\star)^T+ \norm{\tv}_2^2{\bf 1}{\bf 1}^T.
\end{align}

Following \eqref{eq08}, we can estimate the permutation by matching the eigenvectors of $ \Pv\Pv^T-\norm{\tv}_2^2{\bf 1}{\bf 1}^T$ and $\Qv\Qv^T$, complemented by an additional translation estimation procedure. Our algorithm consists of the following steps.
\begin{enumerate}
    \item \textbf{Translation removal}: To isolate the rotational component from the translation vector, we first remove the mean of each point in the data matrices. By multiplying $\frac{1}{m}{\bf 1}{\bf 1}^T$ on both sides of \eqref{eq_irmodel}, we have 
\begin{align}\label{eq_irmodel_avg}
    \frac{1}{m}{\bf 1}{\bf 1}^T\Pv=\frac{1}{m}{\bf 1}{\bf 1}^T\Qv\Rv+{\bf 1}\tv^T.
\end{align}
Subtracting \eqref{eq_irmodel_avg} from  \eqref{eq_irmodel}, we have
\begin{align}\label{eq27}
    \underbrace{\Pv-\frac{1}{m}{\bf 1}{\bf 1}^T\Pv}_{\triangleq \widetilde\Pv}=\Piv^\star \underbrace{\left(\Qv-\frac{1}{m}{\bf 1}{\bf 1}^T\Qv\right)}_{\triangleq \widetilde\Qv}\Rv.
\end{align}
\item \textbf{Permutation estimation}:
The model in \eqref{eq27} now aligns with the one in \eqref{eq01}. We compute $\widehat \Piv$ by using the spectral matching algorithm in \eqref{eq08} with $\Vv$ and $\Uv$ representing the eigenvector matrices of $\widetilde\Qv\widetilde\Qv^T$ and $\widetilde\Pv\widetilde\Pv^T$, respectively.
\item \textbf{Rotation estimation}: Given $\widehat \Piv$, \eqref{eq_ir} becomes the LS problem over the rotation group:
\begin{align}
        \widehat\Rv=\argmin_{\Rv\in SO(3)}\norm{\widetilde\Pv-\widehat\Piv\widetilde\Qv\Rv}_F^2=\argmax_{\Rv\in SO(3)}\tr(\widetilde\Pv^T\widehat\Piv\widetilde\Qv\Rv).
\end{align}
Denote the left and right singular matrices of the matrix $\widetilde\Pv^T\widehat\Piv\widetilde\Qv$ as 
$\Uv_{\text{SVD}}$ and $\Vv_{\text{SVD}}$, respectively. The solution is given by
\begin{align}
    \widehat\Rv=\Vv_{\text{SVD}}\Uv_{\text{SVD}}^T.
\end{align}

\item \textbf{Translation estimation}: The translation is computed by 
\begin{align}
    \hat\tv=  \frac{1}{m}\left(\Pv-  \Qv\widehat\Rv\right)^T{\bf 1}.
\end{align}
\end{enumerate}

\remark[\it Comparisons with alternating optimization]{
The proposed approach, while seemingly analogous to the popular alternating optimization algorithms used in solving \eqref{eq_ir}, diverges fundamentally in its methodology. Alternating optimization typically begins with randomly initialized estimates for $\Qv$ and $\tv$, and iteratively updates $\Piv$, $\Rv$, and $\tv$ until convergence. This method relies heavily on the quality of the initialization and is susceptible to converging to local minima. In contrast, our algorithm is non-iterative and requires no initial guesses. Thus, our solution is simpler and more robust.
}
\section{Numerical Results}\label{sec6}
In this section, we present numerical results to evaluate the performance of the proposed approach.

\subsection{Simulation Results on Shuffled LS}
In this section we study the performance of the shuffled LS solution for the observation model in \eqref{eq01}. As mentioned above, our focus is particularly on the regime in which a substantial number of measurement samples is available, i.e., $t$ is large.

The simulation setup is outlined as follows. Measurements are generated according to the model in \eqref{eq01}. Each $\xv_i$ is i.i.d. drawn from the standard Gaussian distribution $\Norm({\bf 0},\Iv_n)$, $\Av$ is composed of i.i.d. column-normalized Gaussian entries, and the noise $\nv_i$ is drawn from the distribution of $\Norm({\bf 0},\sigma^2\Iv_m)$. The noise power is specified by the signal-to-noise ratio (SNR), defined as $1/\sigma^2$. Additionally, as specified in \eqref{eq01}, the measurements $\{\yv_i\}_{\forall i}$ exhibit a row permutation by $\Piv^\star$. This permutation is implemented by randomly selecting $p_em$ rows for random shuffling, while the remaining $(1-p_e)m$ rows are left unpermuted. Consequently, on average, a fraction of $p_e$ of the rows are permuted, with $\frac{1}{m}\E[\tr(\Piv^\star)]=1-p_e$. The parameter $p_e$ is adjusted to manage the shuffling level in the experiments. For the purposes of emphasizing the large number of available measurements, we set $m=200$, $n=100$, and $t=10^4$ unless stated otherwise.

We evaluate the performance of the solution in \eqref{eq08} by comparing it to the following baseline methods, each distinct in their approach to estimating the permutation:
\begin{itemize}
    \item \textbf{Oracle Bound with Known Permutation}: Assuming $\Piv^\star$ is perfectly known, the problem in \eqref{eq02} simplifies to a classical LS problem. The optimal solution is given by $\widehat\Xv=\Av^\dagger(\Piv^\star)^T\Yv$. This baseline represents the best possible performance for shuffled LS.
     \item \textbf{Method in} \cite[Section 2.3]{SLR_ref2}: Designed primarily for a noiseless measurement model, this method computes the permutation as follows:
     \begin{align}
         \widehat\Piv=\argmax_{\Piv\in\mathcal{P}_m} \norm{\text{Lever}(\Yv)-\Piv\cdot\text{Lever}(\Av)}_2^2,
     \end{align}
     where $\text{Lever}(\Yv)\in\Real^{m\times 1}$ denotes the leverage scores of input matrix $\Yv$, given by the squared norms of the rows of its reduced left singular matrices.
     
    \item \textbf{Method in} \cite{pmlr-v119-zhang20n}: This method estimates the permutation by solving the linear assignment problem in \eqref{eq17}. As demonstrated in Lemma \ref{lemma1}, this method can achieve a diminishing estimation error when the permutation level is limited, i.e., $p_e \leq 0.25$.
    
 \item \textbf{Method in} \cite{Alternating_Mat_SLR}: This method solves the following problem:
        \begin{align}\label{eq34}
         \widehat\Piv=\argmax_{\Piv\in\mathcal{P}_m}\tr\left(\Piv\Av\cdot thres(\Av^T\Yv)\cdot\Yv^T\right),
        \end{align}
        where $thres(\cdot)$ applies a thresholding operation to each column, setting all entries to zero except for the one with the largest magnitude.
\end{itemize}
After estimating the permutation, all the methods proceed to estimate $\hat\Xv$ as per \eqref{eq04}. In our simulations, we adopt the Hungarian method \cite{Hungarian} to solve all the linear assignment problems. For the proposed method, we select the eigenvectors indexed by $\mathcal{A}$ in \eqref{eq10} that correspond to spectral gaps larger than $10^{-3}$ to compute \eqref{eq08}.

We assess the estimation performance using error metrics relevant to both $\Piv^\star$ and $\Xv$. Specifically, the permutation error is quantified by the average error rate, which ranges from $[0,1]$ and is defined as:
\begin{align}
    &\text{Permutation error rate}\triangleq \frac{1 }{m}\sum_{j=1}^m\mathbbm{1}_{\{\exists k: [\Piv^\star]_{j,k}\neq [\widehat\Piv]_{j,k}\}},
\end{align}
where $\mathbbm{1}_{\zeta}$ is the indicator function equal to one if event $\zeta$ is true and zero otherwise.
The estimation error of $\Xv$ is assessed using the normalized mean squared error (NMSE) defined as:
\begin{align}
    &\text{NMSE of }\Xv\triangleq \frac{\norm{\Xv^\star-\widehat\Xv}_F^2}{nt}.
\end{align}
We perform $300$ Monte Carlo trials and report the average over all the trials in the simulations.

\begin{figure}[!t]
	\centering
	\includegraphics[width=3.2 in]{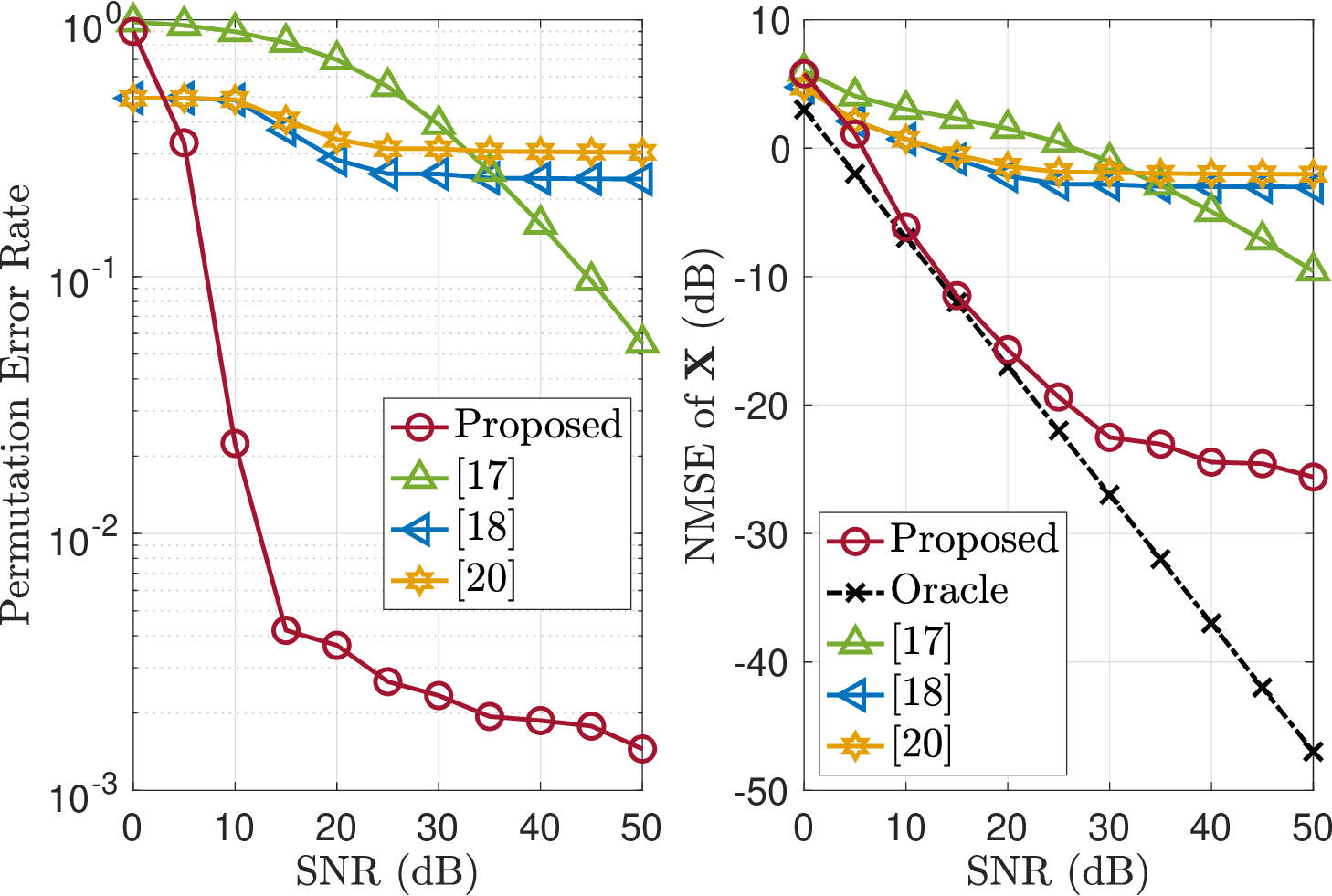}
	\caption{Shuffled LS performance with varying SNR values.}
	\label{fig_sls_snr}
\end{figure}
\begin{figure}[!t]
	\centering
	\includegraphics[width=3.2 in]{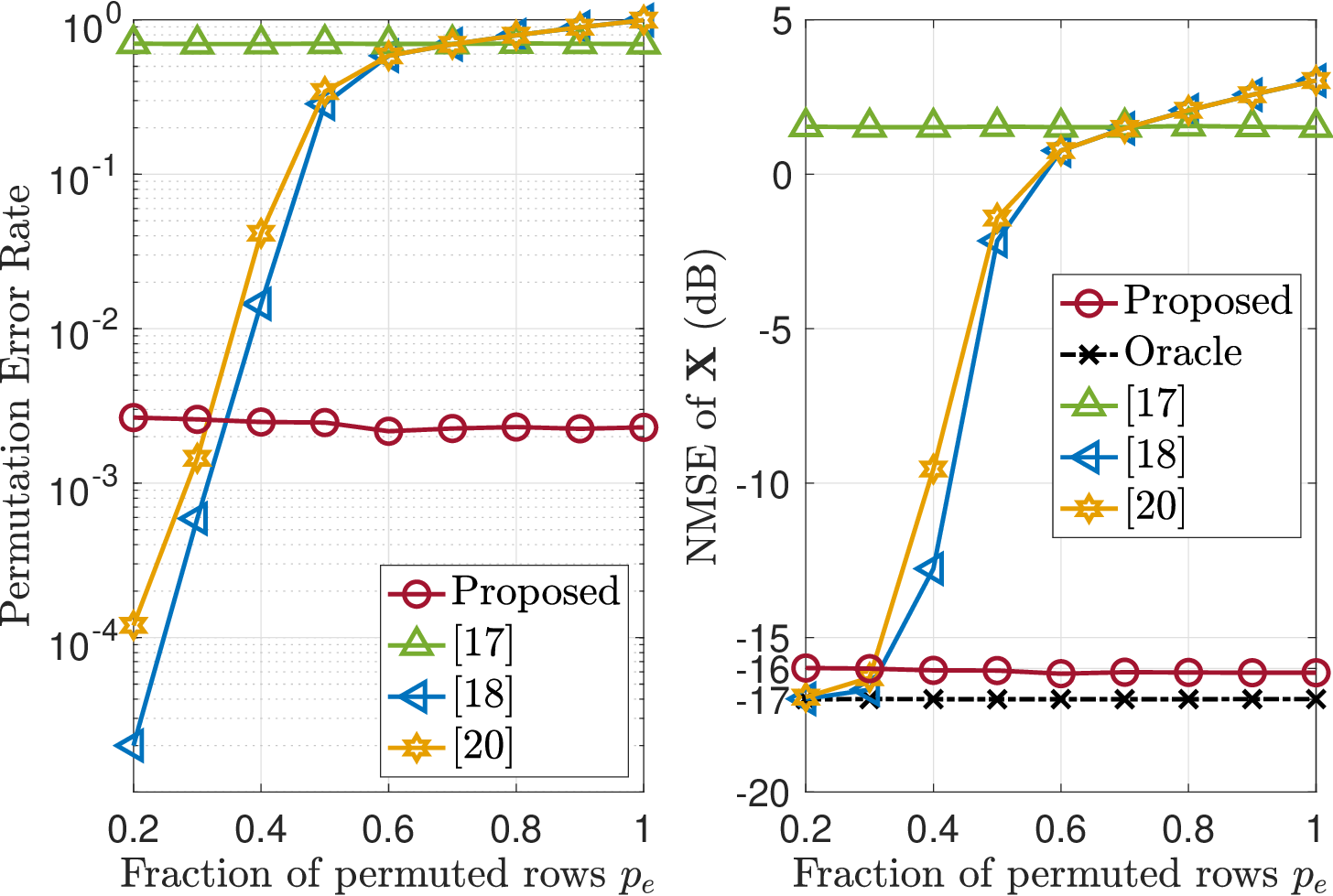}
	\caption{Shuffled LS accuracy versus the fraction of permuted rows $p_e$ in the observations with $\text{SNR}=20$ dB.}
	\label{fig_sls_pe}
\end{figure}

First, we examine the impact of the SNR on estimation accuracy. Considering that the baseline algorithms from \cite{pmlr-v119-zhang20n,Alternating_Mat_SLR} are most effective in scenarios with a limited number of row permutations in the true permutation $\Piv^\star$, we set the parameter $p_e$ in $\Piv^\star$ to 0.5, indicating that $50\%$ of the rows in observations are randomly permuted. Fig. \ref{fig_sls_snr} displays the permutation error and the estimation NMSE as SNR varies. As the SNR increases, the proposed method achieves a significantly smaller permutation error compared to the baselines, thereby enhancing the estimation of the hidden feature matrix $\Xv^\star$. However, in the extremely high SNR regime, the algorithm experiences sporadic permutation errors, leading to an error floor in terms of NMSE of $\Xv$. In contrast, the baseline algorithms \cite{pmlr-v119-zhang20n,Alternating_Mat_SLR} exhibit high permutation errors throughout. The method in \cite{SLR_ref2} is designed primarily for noiseless setups and works only in the extremely high SNR regime.

\begin{figure}[!t]
	\centering
	\includegraphics[width=3.2 in]{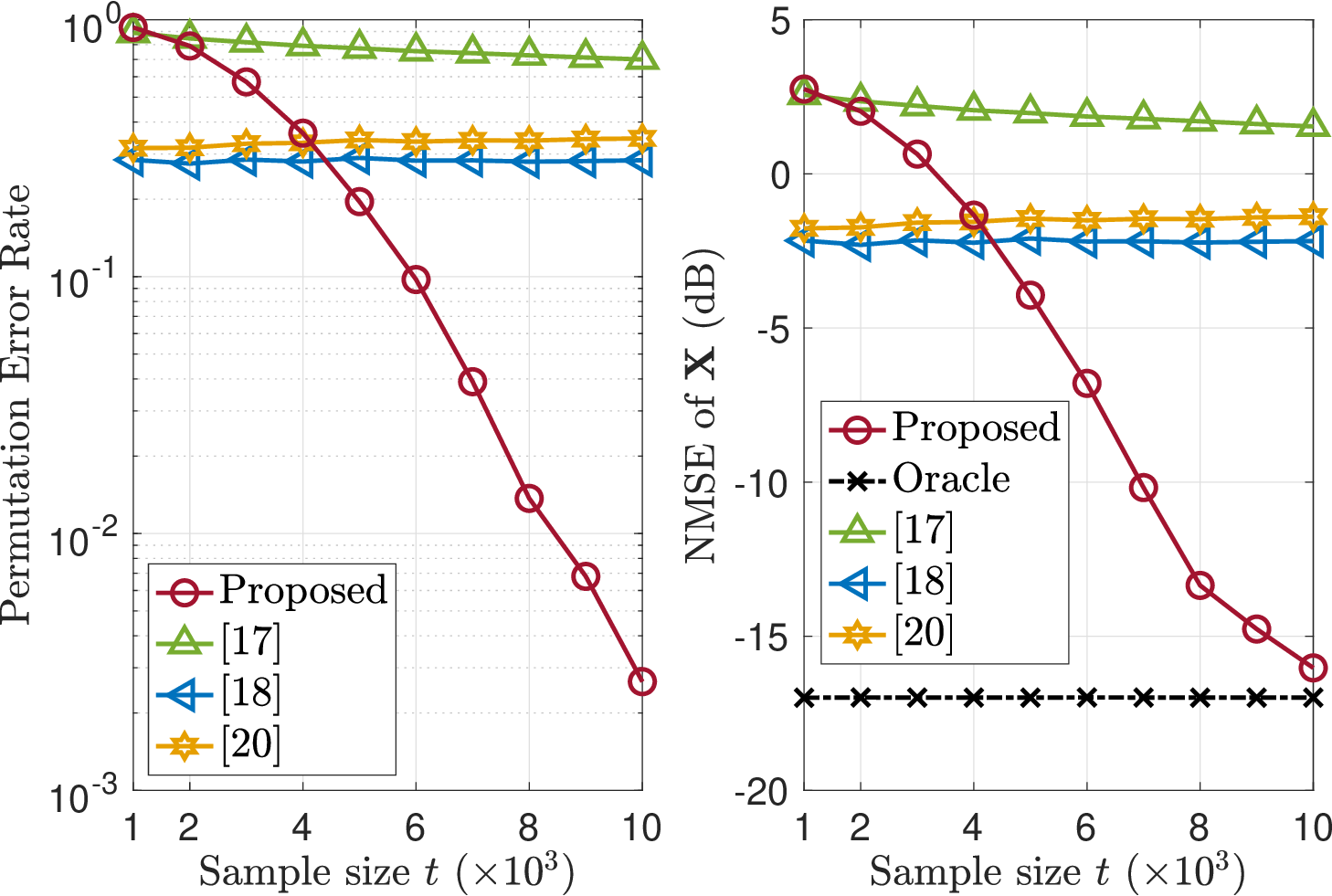}
	\caption{Shuffled LS accuracy versus the sample size $t$ with $p_e=0.5$ and $\text{SNR}=20$ dB.}
	\label{fig_sls_t}
\end{figure}
\begin{figure}[!t]
	\centering
	\includegraphics[width=3.2 in]{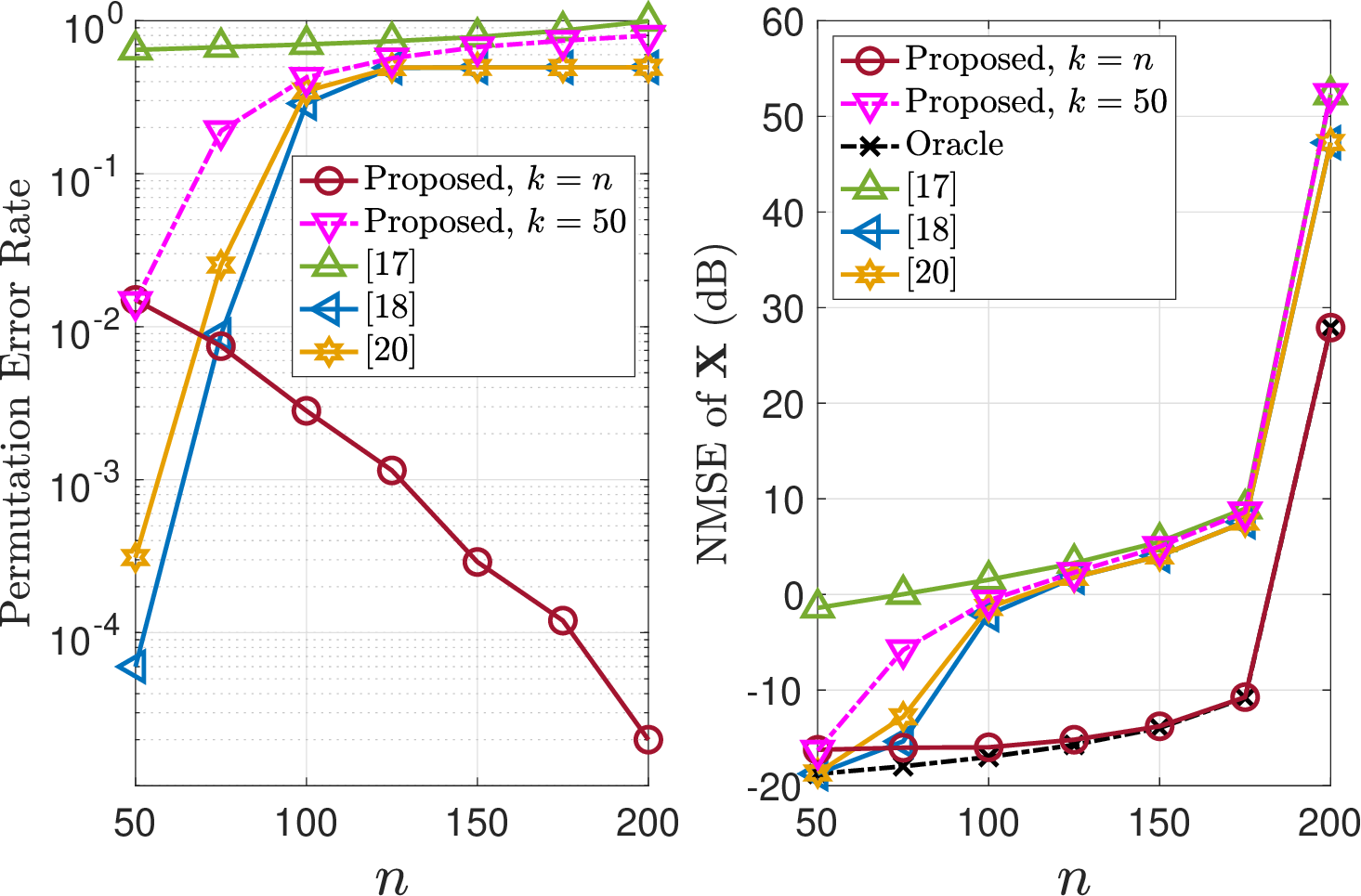}
	\caption{Shuffled LS accuracy versus $n$ with $m=200$, $\text{SNR}=20$ dB, and $p_e=0.5$.}
	\label{fig_sls_n}
\end{figure}

Next, we set the SNR to $20$ dB and vary the fraction of permuted rows in the observations by adjusting the value of $p_e$, as illustrated in Fig. \ref{fig_sls_pe}. As detailed in Section \ref{sec4}, our method does not require a predetermined number of permutations and is capable of identifying any random permutations. The estimation accuracy of our proposed algorithm remains robust to changes in $p_e$. In contrast, the accuracy of baseline algorithms from \cite{Alternating_Mat_SLR, pmlr-v119-zhang20n} significantly worsens as $p_e$ increases. These baseline methods outperform our proposed algorithm when $p_e < 0.3$, yielding gains of $1$ dB to $2$ dB in the estimation NMSE of $\Xv$. However, when $p_e > 0.3$, our method demonstrates markedly superior estimation accuracy. Note that in many real-world scenarios, prior knowledge of the true permutation level is often unavailable. The results in Fig. \ref{fig_sls_pe} underscore the robustness and practical relevance of our method.

Fig. \ref{fig_sls_t} examines the impact of the number of measurements $t$, with fixed values $m=200$ and $n=100$. As $t$ increases, the estimation error of the proposed algorithm decreases significantly and approaches the oracle bound for $t \geq 10^4$. This observation is consistent with the analysis presented in Theorems \ref{theorem1}--\ref{theorem3}, highlighting the necessity of a sufficiently large $t$ for our method to achieve diminishing error. Conversely, as indicated in Lemma \ref{lemma1}, the permutation estimation accuracy of the baseline algorithm \cite{pmlr-v119-zhang20n} does not improve with an increase in $t$ once $t \geq m$.

Finally, with fixed parameters $m=200$ and $t=10^4$, we explore the impact of varying $n$ within the range $[50,200]$. 
According to Theorem \ref{theorem3}, the number of selected eigenvectors $k$ should scale at least linearly with $n$. We consider two choices for $k$ in the simulations: setting $k=n$ and $k=50$, depicted by the dashed magenta and solid red curves, respectively, in Fig. \ref{fig_sls_n}. These results confirm the theoretical findings in Theorem \ref{theorem3}, illustrating the importance of the condition $k \gtrsim n$ for achieving diminishing estimation error in our shuffled LS method. This condition is critical for ensuring that our method surpasses the performance of the baseline algorithms.

\subsection{Simulation Results on Shuffled LASSO}
We examine the performance of the shuffled LASSO solutions in \eqref{eq03}. As in the shuffled LS setup, measurements are generated according to \eqref{eq01} using a column-normalized sensing matrix $\Av$ and an i.i.d. noise matrix $\Nv$, with entries distributed according to $\Norm(0,\sigma^2)$. To incorporate a sparse feature matrix, we model the entries of $\Xv$ using an i.i.d. Bernoulli-Gaussian distribution: 
\begin{align}
    p(x_{i,j})=\left(1-\frac{s}{n}\right)\mathbbm{1}_{\{x_{i,j}=0\}}+\frac{s}{n}\Norm(0,1),\forall i,j,
\end{align}
where $s$ is a hyperparameter that controls the sparsity level, ensuring that $\E[\norm{\xv_i}_0] = s$ for all $i$. Additionally, the true permutation matrix is uniformly drawn from the set of permutation matrices $\mathcal{P}_m$.

\begin{figure}[!t]
	\centering
	\includegraphics[width=3.2 in]{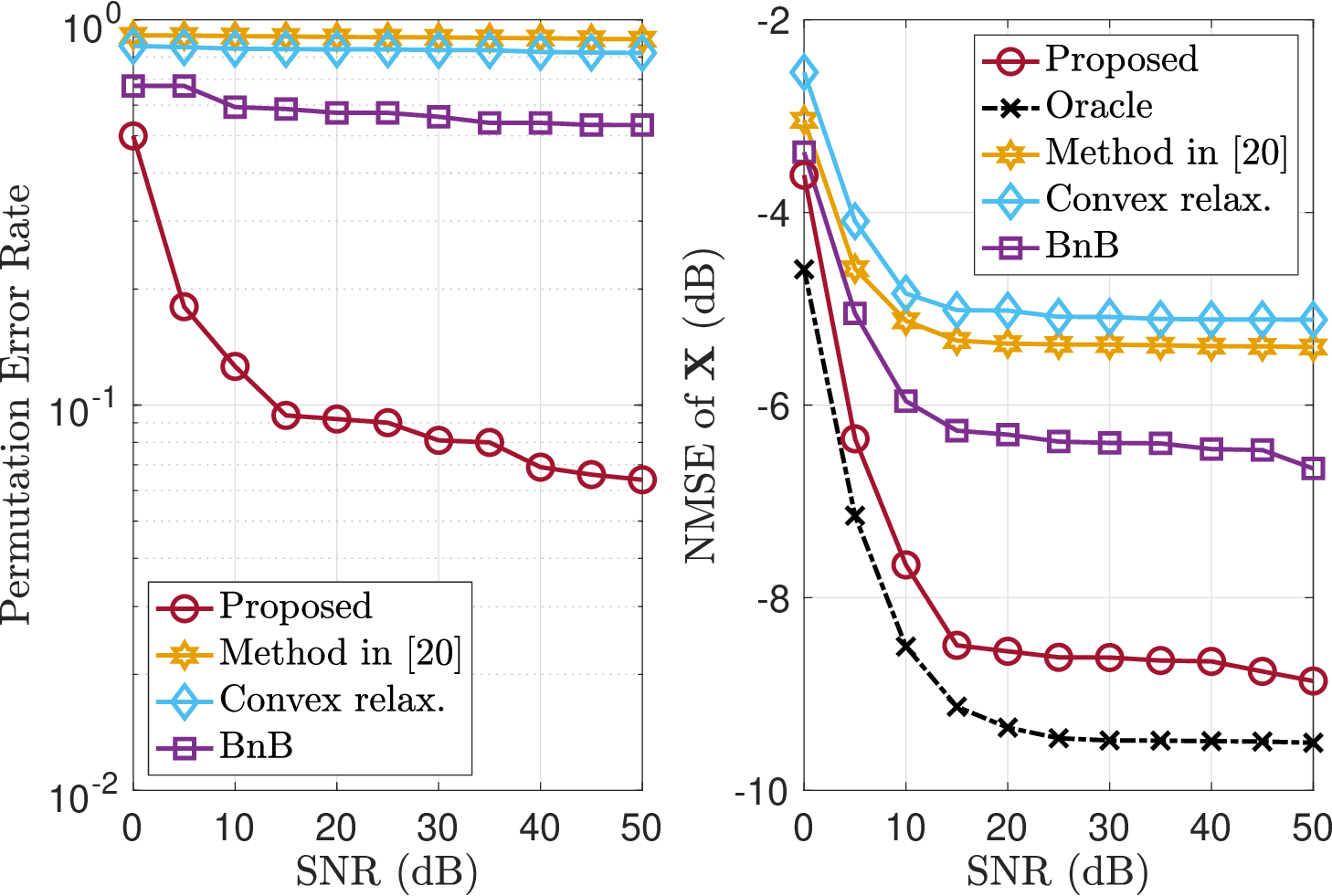}
	\caption{Shuffled LASSO performance versus SNR.}
	\label{fig_lasso_snr}
\end{figure}

We evaluate the performance of the proposed solution in \eqref{eq08} and \eqref{eq04b} by comparing it to the following baselines:
\begin{itemize}
    \item \textbf{Oracle Bound}: Assuming $\Piv^\star$ is known, the shuffled LASSO problem in \eqref{eq03} simplifies to a classical LASSO problem as shown in \eqref{eq04b} with $\widehat\Piv=\Piv^\star$. 
    \item \textbf{Method in} \cite{Alternating_Mat_SLR}: This baseline computes the permutation $\widehat\Piv$ by using \eqref{eq34} and subsequently applies the LASSO solution in \eqref{eq04b}.
    \item \textbf{Convex Relaxation} \cite[Algorithm 1]{SLR_CS3}: This method addresses the problem through an $\ell_1$-minimization approach:\footnote{Originally, the approach in \cite{SLR_CS3} was designed for a noiseless environment. We have adapted it to accommodate noisy observations by modifying the linear constraint to an $\ell_2$-norm constraint, which aligns with common practices in compressed sensing \cite{10.5555/2526243}.}
    \begin{align}
        \min_{\Xv,\Piv}&~~ \norm{\Xv}_{1,1}\nonumber\\
        \text{s.t. }& \norm{\Yv-\Piv\Av\Xv}_F\leq \eta,\Piv\in\mathcal{P}_m,\label{eq37}
    \end{align}
    where $\eta$ is a predefined parameter. Given the combinatorial nature of $\Piv$, this non-convex problem poses significant challenges. To tackle this challenge, Ref. \cite{SLR_CS3} relaxes the feasible set for $\Piv$ from the set of permutation matrices to that of doubly-stochastic matrices, transforming the problem into a convex one.
After solving the resultant convex programming, the solution for the doubly-stochastic matrix is projected back onto the set of permutation matrices $\mathcal{P}_m$. Subsequently, $\Xv$ is estimated by addressing the $\ell_2$-constrained $\ell_1$-minimization problem using the estimated permutation matrix.
    \item \textbf{Branch-and-Bound (BnB) Algorithm} \cite[Algorithm 3]{SLR_CS3}: This method tackles the non-convex problem described in \eqref{eq37} by employing the BnB algorithm. This approach iteratively divides the permutation set into subsets and searches for the optimal solution within them.
\end{itemize}
Unless otherwise specified, in our simulations we fix $n=30$, $m=10$, $t=1000$, and $s=6$. This choice results in an underdetermined system, even when permutations are known, thereby underscoring the critical role of exploiting feature sparsity. { According to Theorem \ref{theorem4}, a sufficient condition for achieving diminishing estimation error requires the penalty coefficient in \eqref{eq03} to satisfy $\rho \gtrsim \sigma \sqrt{\frac{\ln n}{m}}$. Motivated by this, we adjust the penalty coefficient as $\rho = 20 \sigma \sqrt{\frac{\ln n}{m}}$.\footnote{ Alternatively, when knowledge of system parameters such as $\sigma$ is unavailable, the penalty parameter can be selected using conventional parameter tuning methods, such as grid search with cross-validation or feature selection via Akaike Information Criterion (AIC) or the Bayesian Information Criterion (BIC).} For the baseline algorithms in \cite{SLR_CS3}, we tested various values for the hyperparameter $\eta$ in \eqref{eq37} and found that $\eta = 10$ generally yields the best results across most setups.} All the convex programming problems are solved using CVX \cite{cvx}.

\begin{figure}[!t]
	\centering
	\includegraphics[width=3.2 in]{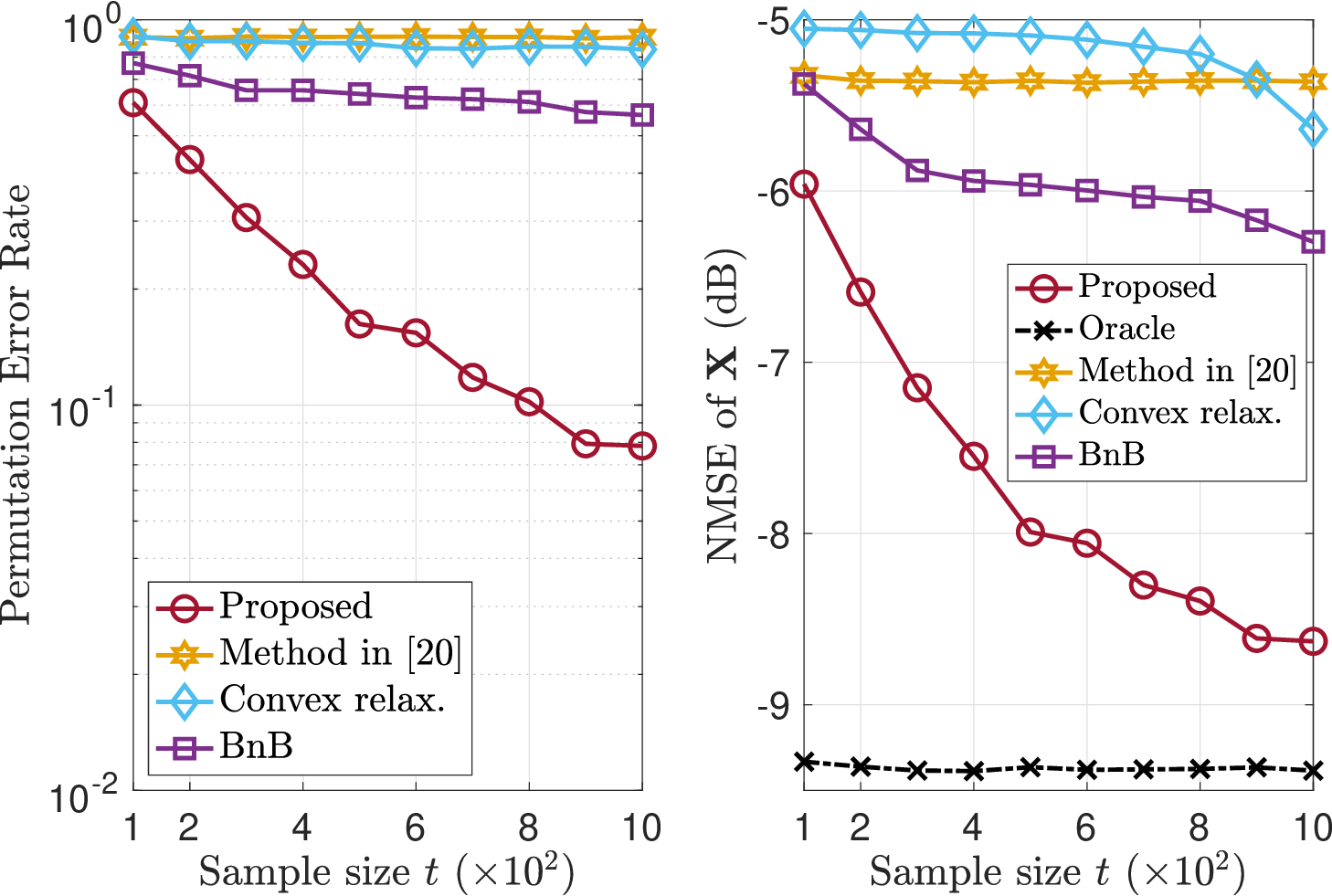}
	\caption{Shuffled LASSO performance versus sample size $t$.}
	\label{fig_lasso_t}
\end{figure}

\begin{figure}[!t]
	\centering
	\includegraphics[width=3.2 in]{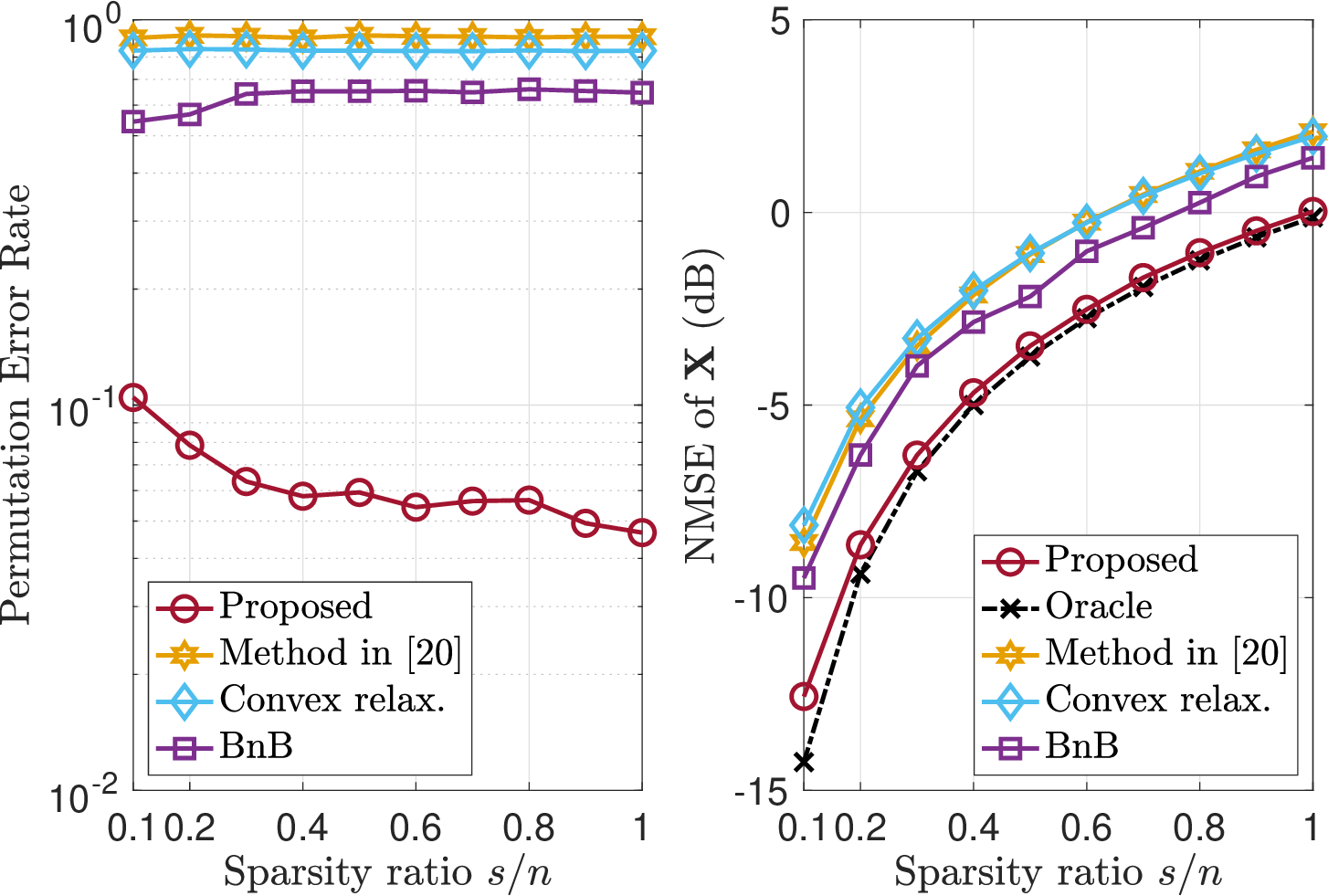}
	\caption{Shuffled LASSO accuracy versus the sparsity ratio.}
	\label{fig_lasso_s}
\end{figure}
Fig. \ref{fig_lasso_snr} displays the estimation accuracy of the shuffled LASSO across varying SNR values. The proposed approach outperforms all the baseline algorithms, primarily due to its more precise estimation of permutations. In Fig. \ref{fig_lasso_t}, we maintain the SNR at $20$ dB and vary the sample size $t$. For the proposed method and the $\ell_1$-minimization approaches from \cite{SLR_CS3}, increasing $t$ leads to more accurate permutation estimation and, consequently, a better estimation of the underlying feature matrix $\Xv$. { Conversely, the permutation estimation accuracy of the baseline algorithm in \cite{Alternating_Mat_SLR} is insensitive to the changes in $t$ once $t \geq m$. This observation is consistent with the results in Fig. \ref{fig_sls_t}.}

\begin{figure}[!t]
	\centering
	\includegraphics[width=3 in]{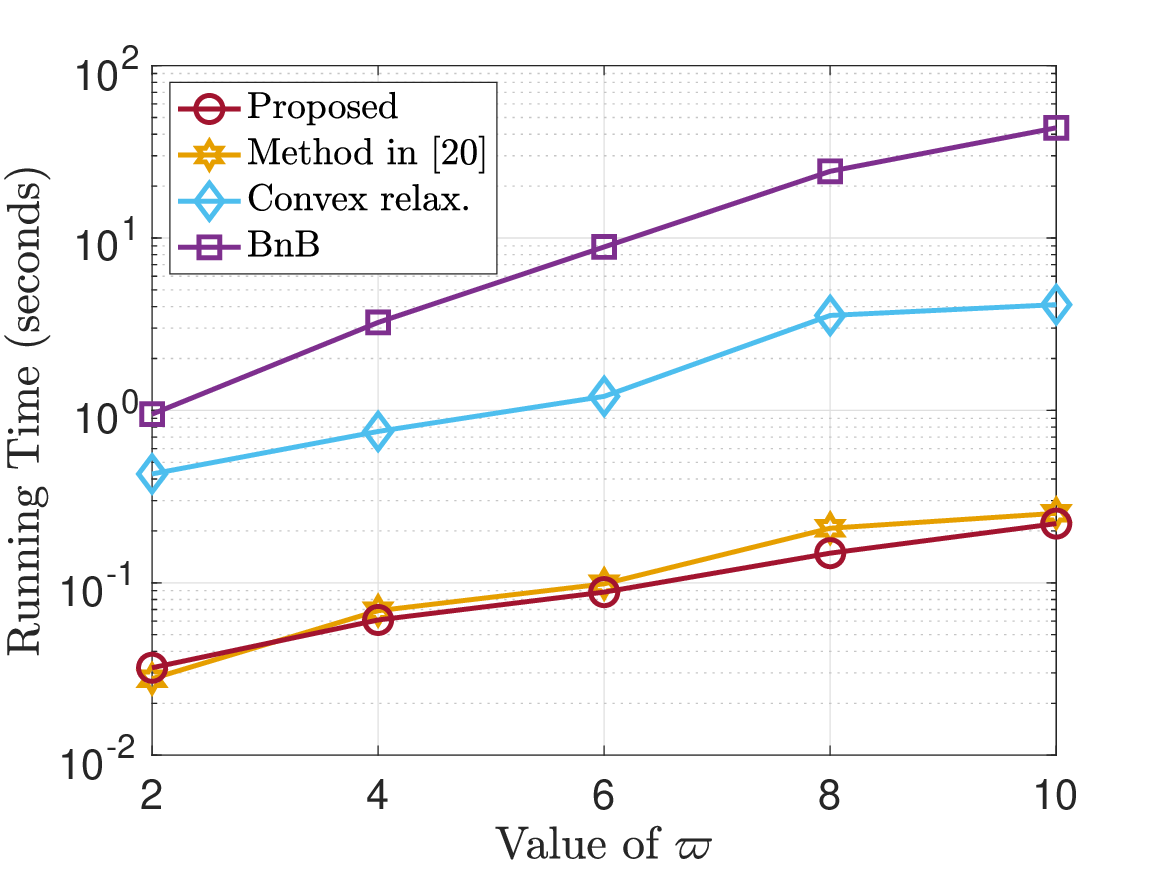}
	\caption{Average running time (in seconds) of different matching algorithms for the shuffled LASSO problem versus the problem size, where the system size parameters are given by $m=\varpi$, $n=3\varpi$ and $t=100\varpi$.}
	\label{fig_lasso_size}
\end{figure}

As stated in Theorem \ref{theorem4}, with a fixed system size, the estimation MSE is bounded by a term proportional to the sparsity of the feature vectors. To explore the impact of feature sparsity, Fig. \ref{fig_lasso_s} illustrates performance relative to the average sparsity level, represented by ${s}/{n}$. We observe that the permutation estimation accuracy of our method is insensitive to the features' sparsity. This is because our spectral matching approach does not rely on feature sparsity for accurate permutation estimation and is thus suitable for both sparse and dense feature sets. Moreover, the oracle bound, assuming $\Piv^\star$ is known, shows that the estimation NMSE of the feature matrix increases with denser features. Remarkably, the proposed method approaches this oracle bound in estimation accuracy, underscoring the effectiveness of our spectral matching technique.

{ Finally, we conduct simulations to assess the running time of the proposed method, validating the computational complexity analysis. We simulate the shuffled LASSO problem with varying problem sizes by proportionally scaling the size parameters $m$, $n$, and $t$. Specifically, we use a system size parameter $\varpi$ to adjust these parameters as follows: $m = \varpi$, $n = 3\varpi$, and $t = 100\varpi$, where a larger $\varpi$ corresponds to a higher-dimensional matching system. The other parameters remain the same as those in Fig. \ref{fig_lasso_s} with $s/n=0.2$.  Fig. \ref{fig_lasso_size} plots the average running time (in seconds) of different matching algorithms versus the value of $\varpi$, averaged over $50$ Monte Carlo trials. The simulations were conducted using MATLAB on a macOS computer with a 2.6 GHz CPU and 16 GB of memory.  The results show that the baseline algorithms in \cite{SLR_CS3}, particularly the BnB algorithm, exhibits significantly higher time costs in large-scale systems. In contrast, the proposed method and the baseline algorithm in \cite{Alternating_Mat_SLR} requires less running time. This observation is consistent with the complexity analysis presented in Section \ref{sec3}.}
\begin{figure*}[!t]

	\centering
	\includegraphics[width=7 in]{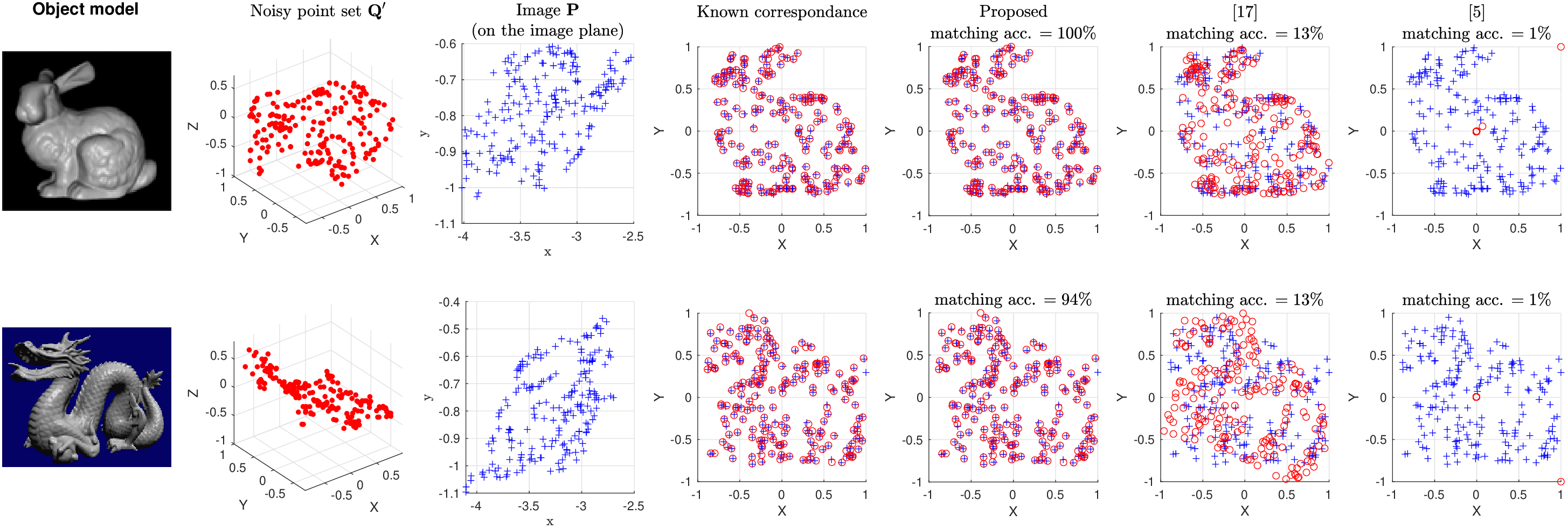}
	\caption{Image registration for the `Stanford Bunny' (Top) and `Dragon' (Bottom) models. From left to right: 1) The 3D scanned object models; 2) The observed noisy 3D point set $\Qv^\prime$; 3) The 2D image point set $\Pv$; 4)--7) Registration results for the oracle bound with known correspondence, the proposed method, and the baseline algorithms, visualized on the $XY$-plane of the original object coordinates. The matching accuracy measures the fraction of correctly corresponding points.
 } 
 \hrulefill
	\label{fig_3dir}
\end{figure*}

\begin{figure}[!t]
	\centering
	\includegraphics[width=3.2 in]{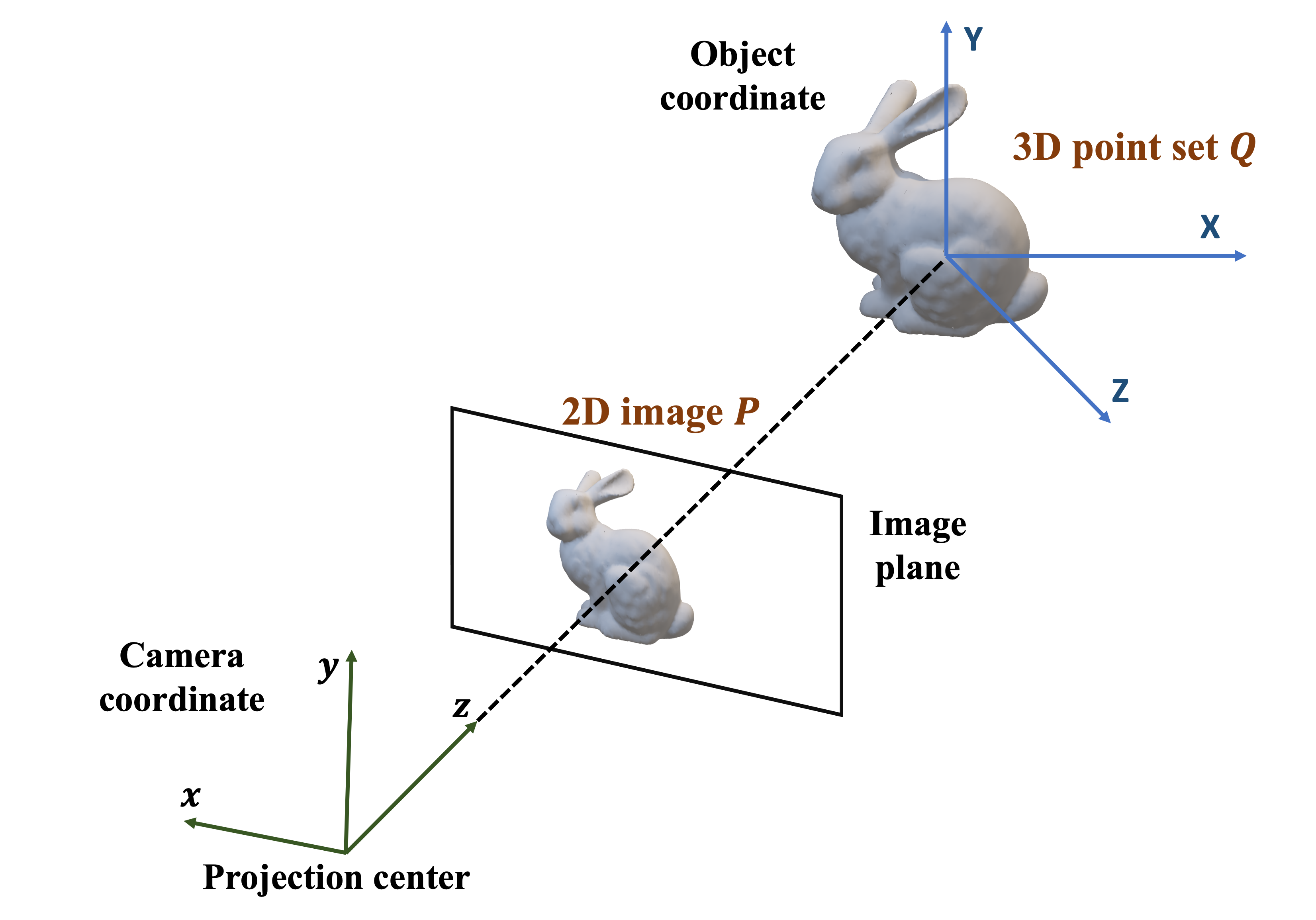}
	\caption{Illustration of the image registration setup.}
	\label{fig_3Dillustration}
\end{figure}

\subsection{Image Registration Results}\label{sec6c}
In this section, we apply the proposed SLR algorithm to address the simultaneous pose and correspondence problem within the context of 3D-2D image registration. We utilize the real-world dataset from the Stanford 3D Scanning Repository \cite{Stanford3D}, which includes 3D models of scanned objects.

The 3D-2D registration challenge is structured as per the methodology described in \cite{softposit}, with the mathematical relationship expressed in \eqref{eq_irmodel}. For each scanning model, a (noiseless) 3D point set $\Qv$ is generated by randomly sampling $200$ points from the model data. The Cartesian coordinates of these points are normalized within the range $[0,1]^3$. As shown in Fig. \ref{fig_3Dillustration}, the corresponding 2D image set $\Pv$ is produced by projecting $\Qv$ with the projection center located at $[10,10,0]$, which translates to a vector $\tv=[10,10,0]$ in \eqref{eq_irmodel}. Originating from this projection center, a camera coordinate system is established with randomly drawn orthogonal $xyz$-axes, corresponding to the unknown rotation $\Rv$. The image plane is aligned along the resultant $z$-axis. Due to the loss of point correspondence during projection, the set $\Pv$ is organized by aligning the homogeneous coordinates of the projected points from the top-left to the bottom-right in the image plane.

Suppose we observe $\Pv$ and a noisy version of the original 3D point set, denoted by $\Qv^\prime$, whose coordinates are perturbed by i.i.d. zero-mean Gaussian noise with a variance $\sigma^2=10^{-4}$. Our objective is to register $\Pv$ to this noisy set $\Qv^\prime$. This involves accurately estimating the rotation $\Rv$, the translation $\tv$, and the correspondence $\Piv^\star$ using the shuffled LS framework outlined in \eqref{eq_ir}. We compare the registration performance of our proposed algorithm, as detailed in Section \ref{sec5}, against the baseline algorithms presented in \cite{softposit} and \cite{SLR_ref2}, as well as an oracle baseline assuming the point correspondence is perfectly known.
Fig. \ref{fig_3dir} plots the registration results for the `Stanford bunny' and `dragon' models. To aid in the comparison, we visualize the estimated point correspondence between the reconstructed 3D point set $\widehat\Qv = \widehat\Piv \Pv \widehat\Rv + {\bf 1} \hat\tv^T$ and the noisy observation $\Qv^\prime$ on the $XY$-plane of the original object coordinates. As discussed in Section \ref{sec5}, this image registration problem does not meet the condition $t \gg m$, which is essential for accurate spectral matching for i.i.d. features as described in Section \ref{sec3}. Nonetheless, the orthogonality of the feature matrix, i.e., the rotation matrix in this context, supports correct covariance matrix estimation even when $t = 3$. As illustrated in Fig. \ref{fig_3dir}, our proposed method significantly outperforms the baseline algorithms, achieving more accurate correspondence estimation and a near-optimal point set reconstruction after registration.

\section{Conclusions}\label{sec7}
In this work, we explored the SLR framework designed to estimate hidden features from linear transformations under unknown permutations. We delved into both shuffled LS and shuffled LASSO problems, focusing particularly on scenarios involving a large number of measurement samples. We proposed a spectral matching approach for permutation estimation and further extended this method to address simultaneous pose and correspondence estimation problems in image registration. Our theoretical analysis provided insights into the accuracy of permutation estimation and the feature estimation MSEs for the resulting shuffled LS and LASSO solutions. In particular, we established that the estimation error of our algorithm converges to an asymptotically optimal error rate when the sample size is sufficiently large. Numerical results on both synthetic datasets and real-world image registration applications corroborate our theoretical findings and demonstrate that our method outperforms existing state-of-the-art algorithms.

\appendices
\section{Definitions and Notations}\label{app_def}
In this section, we summarize the notations used in our proofs. 
For any permutation matrix $\Piv\in\mathcal{P}_m$, we define a corresponding index mapping function $\pi(\cdot):\{1,\cdots,m\}\to \{1,\cdots,m\}$ such that $\pi(k)=l$ if and only if $[\Piv]_{k,l}=1$. We shall use $\Piv$ and $\pi(\cdot)$ to denote the permutation interchangeably. We denote the index mapping corresponding to $\Piv^\star$ and $\widehat \Piv$ by $\pi^\star$ and $\hat \pi$, respectively. 

For any given $|\mathcal{A}|$, we denote $k=|\mathcal{A}|$. Recall that the eigenvalues of $\Lambdav$ and $\widehat\Lambdav$ in \eqref{eq06} and \eqref{eq07} are aligned in the descending order, \ie $\lambda_1>\cdots\lambda_m$ and $\hat\lambda_1\geq \cdots\geq \hat \lambda_m$. We use $\delta_i=\max\{\lambda_{i+1}-\lambda_i,\lambda_i-\lambda_{i-1}\}$ to denote the $i$-th spectral gap of $\Lambdav$ with $\lambda_0\triangleq \infty$ and $\lambda_{m+1}\triangleq 0$. 
We denote the minimum spectral gap in $\mathcal{A}$ by $\delta(\mathcal{A})=\min_{i\in\mathcal{A}} \delta_i$. For simplicity, when the context is clear, we omit the argument $\mathcal{A}$ and simply denote $\delta(\mathcal{A})$ as $\delta$.


Finally, we denote the cost matrices in \eqref{eq08} and \eqref{eq_noiselesspro} by $\widehat\Dv=abs\left(\Vv_\mathcal{A})abs(\Uv_\mathcal{A}^T\right)$ and $\Dv=\Vv_\mathcal{A}(\Piv^\star\Vv_\mathcal{A})^T$, respectively. As a consequence, the problems in \eqref{eq08} and \eqref{eq_noiselesspro} can be recast as
\begin{align}\label{eq_appdef1}
     \widehat\Piv\in\argmax_{\Piv\in\mathcal{P}_m}\tr\left(\widehat \Dv\Piv \right),
\end{align}
and 
\begin{align}\label{eq_appdef2}
     \Piv^\star\in\argmax_{\Piv\in\mathcal{P}_m}\tr\left(\Dv\Piv \right).
\end{align}

\section{Proof of Theorem \ref{theorem1}}\label{appa}
Using the definition in \eqref{eq_appdef2}, the bound to prove in \eqref{eq14} is equivalent to
\begin{align}
    1-\frac{1}{k}\tr(\Dv\widehat\Piv)\lesssim \frac{m}{\delta}\sqrt{\frac{m(\ln m +\epsilon)}{kt}}.
\end{align}

First, using the definition in \eqref{eq_appdef1}, we have
\begin{align}\label{eq_appa11}
&\tr\left( \Dv\widehat\Piv \right)=\tr\left( \widehat\Dv\widehat\Piv \right)+\tr((\Dv-\widehat\Dv)\widehat\Piv )\nonumber\\
\geq& \tr\left( \widehat\Dv\Piv^\star \right)-m \norm{\Dv-\widehat\Dv}_{\max},
\end{align}
where $\norm{\Dv-\widehat\Dv}_{\max}=\max_{i,j}[\Dv-\widehat\Dv]_{i,j}$ denotes the max-norm. Applying the triangle inequality and the matrix norm inequality, 
\begin{align}
    &\norm{\Dv-\widehat\Dv}_{\max}\nonumber\\
    =&\norm{\Vv_\mathcal{A} \Vv_\mathcal{A}^T (\Piv^\star)^T-abs(\Vv_\mathcal{A})abs(\Uv_\mathcal{A}^T)}_{\max}\nonumber\\
    \leq& \norm{\Vv_\mathcal{A} \Vv_\mathcal{A}^T (\Piv^\star)^T-abs(\Vv_\mathcal{A}\Uv_\mathcal{A}^T)}_{\max}\nonumber\\
    \leq& \norm{\Vv_\mathcal{A}(\Piv^\star\Vv_\mathcal{A}-\Uv_\mathcal{A})^T}_2\leq \norm{\Piv^\star\Vv_\mathcal{A}-\Uv_\mathcal{A}}_F,\label{eq_appa12}
\end{align}
where the last inequality follows from $\norm{\Vv_{\mathcal{A}}}_2=1$.
Hence, 
\begin{align}
    &k-\tr(\Dv\widehat\Piv)\overset{\eqref{eq_appa12}}{\leq} k-\sum_{j=1}^m [\widehat\Dv]_{j,\pi^\star(j)}+m\norm{\Piv^\star\Vv_\mathcal{A}-\Uv_\mathcal{A}}_F\nonumber\\
    =&k-\tr(abs(\Vv_\mathcal{A})abs(\Uv_\mathcal{A}^T)\Piv^\star)+m\norm{\Piv^\star\Vv_\mathcal{A}-\Uv_\mathcal{A}}_F\nonumber\\
    \leq&k-\tr(abs(\Vv_\mathcal{A}\Uv_\mathcal{A}^T)\Piv^\star)+m\norm{\Piv^\star\Vv_\mathcal{A}-\Uv_\mathcal{A}}_F\nonumber\\
    \leq& k-\tr((\Vv_\mathcal{A}\Uv_\mathcal{A}^T\Piv^\star)^T(\Vv_\mathcal{A}\Uv_\mathcal{A}^T\Piv^\star))+m\norm{\Piv^\star\Vv_\mathcal{A}-\Uv_\mathcal{A}}_F\nonumber\\
    =&\frac{1}{2}\norm{\Piv^\star\Vv_\mathcal{A}-\Uv_\mathcal{A}}_F^2+m\norm{\Piv^\star\Vv_\mathcal{A}-\Uv_\mathcal{A}}_F,\label{eq_appa14}
\end{align}
where the last inequality follows from that each entry in $abs(\Vv_\mathcal{A}\Uv_\mathcal{A}^T)\Piv^\star$ lies in the range of $[0,1]$. 

Applying the inequality in \cite[Eq. (54)]{BlindGM} and the Davis-Kahan theorem in \cite[Theorem 2]{DK_theorem2}, it follows that 
\begin{align}\label{eq_appa08}
    \norm{\Piv^\star\Vv_\mathcal{A}-\Uv_\mathcal{A}}_F\leq \frac{2\sqrt{2k}}{\delta}\norm{\Cv_Y-\widehat\Cv_Y}_2,
\end{align}
where $\delta$ is the minimum spectral gap of the principal components defined in Appendix \ref{app_def}.

In addition, the results in \cite[Theorem 5.6.1]{HDimPro} and \cite[Exercise 5.6.4]{HDimPro} lead to a concentration bound on $\norm{\Cv_Y-\widehat\Cv_Y}_2$. Specifically, when $\norm{\yv_i}_2$ is bounded almost surely, for any $\epsilon>0$, it follows with probability greater than $1-2e^{-\epsilon}$ that 
\begin{align}
    \norm{\Cv_Y-\widehat\Cv_Y}_2\leq c_0\left(\sqrt{\frac{m(\ln m+\epsilon)}{t}}+\frac{m(\ln m + \epsilon)}{t}\right),
\end{align}
where $c_0=2(\lambda_1+\norm{C_N}_2)\sup \norm{\yv_i}^2/\sqrt{\E[\norm{\yv_i}^2_2]}$ depends on the upper bound of the measurement norm. For $t\geq m(\ln m +\epsilon)$, it follows with probability greater than $1-2e^{-\epsilon}$ that 
\begin{align}\label{eq_appa13}
    \norm{\Cv_Y-\widehat\Cv_Y}_2\leq 2c_0\sqrt{\frac{m (\ln m + \epsilon)}{t}}. 
\end{align}
Combining \eqref{eq_appa14}--\eqref{eq_appa13}, for any $t\geq m(\ln m+\epsilon)/\min\{1,\delta^2\}$, it follows with probability greater than $1-2e^{-\epsilon}$ that 
\begin{align}
    &1-\frac{1}{k}\tr(\Dv\widehat\Piv)\nonumber\\
    \leq& \frac{16c_0^2}{\delta^2}\frac{m(\ln m+\epsilon)}{t}+\frac{4\sqrt{2}c_0}{\delta}m\sqrt{\frac{m(\ln m+\epsilon)}{kt}}\nonumber\\
    \overset{(a)}{\leq} & \frac{16c_0^2}{\delta}\sqrt{\frac{m(\ln m+\epsilon)}{t}}+\frac{4\sqrt{2}c_0}{\delta}\frac{m}{\sqrt{k}}\sqrt{\frac{m(\ln m+\epsilon)}{t}}\nonumber\\
    \lesssim& \frac{m}{\delta}\sqrt{\frac{m(\ln m +\epsilon)}{kt}},
\end{align}
where $(a)$ follows from $\frac{m(\ln m+\epsilon)}{\delta^2 t}\leq 1$.

\section{Proof of Theorem \ref{theorem2}} \label{appb}
When $\xv_i$ exhibits a rotational invariant distribution, it is known that the eigenmatrix $\Vv$ in \eqref{eq06} is distributed according to Haar measure on the $m$-dimensional orthogonal group, i.e., uniformly distributed over the orthogonal matrix group; see, e.g., \cite{Haar_othogonal}. As a consequence, it follows that the entries of $\Vv$ are identically distributed \cite[Lemma 2.1]{Haar_othogonal} and the result in \cite[Lemma 6.1]{Haar_othogonal} provides the following concentration bound:
\lemma{\label{lemma_appc}Let $\Vv_{\mathcal{A}}\in\Real^{m\times k}$ be a matrix, consisting of $k=|\mathcal{A}|$ columns of a matrix distributed according to Haar measure. Let $\{\av_i\in\Real^m\}_{i=1}^n$ and $\{\bv_i\in\Real^m\}_{i=1}^n$ be two  collections of arbitrary $m$-dimensional vectors. Let $u >0$ be a fixed constant. With probability at least $1-Ck^{2-C^\prime \frac{ku^2}{\ln (2n)}}$, the following holds for $\forall 1\leq i,j\leq n$:
\begin{align}
    &\frac{m}{k}\norm{\Vv_{\mathcal{A}}^T\av_i-\Vv_{\mathcal{A}}^T\bv_i}_2^2\geq (1-u)\norm{\av_i-\bv_i}_2^2,\nonumber
\end{align}
where $C$ and $C^\prime$ are some absolute constants. 
}
\begin{IEEEproof}
    See \cite[Lemma 6.1]{Haar_othogonal}.
\end{IEEEproof}

Let us consider the case where $\widehat \Piv\neq \Piv^\star$ holds. Define $\mathcal{T}\triangleq\{1\leq i\leq m:\hat \pi(i)\neq \pi^\star(i)\}$ to be the misaligned index set. Clearly, in this case of $\widehat \Piv\neq \Piv^\star$, we have $|\mathcal{T}|\geq 1$. Moreover, it follows that
\begin{align}\label{eq_appa02}
	&\tr(\widehat\Dv\widehat\Piv)- \tr(\widehat\Dv\Piv^\star)\geq 0,
\end{align}
where $\widehat\Dv$ is defined in \eqref{eq_appdef1}.
On the other hand, we note that
\begin{align}
  &  \tr(\widehat\Dv\widehat\Piv)- \tr(\widehat\Dv\Piv^\star)\nonumber\\
    =&\tr((\Dv-\widehat\Dv)(\Piv^\star-\widehat\Piv))-\tr(\Dv\Piv^\star)+\tr(\Dv\widehat\Piv)\nonumber\\
    \overset{(a)}{\leq }& 2|\mathcal{T}|\norm{\Dv-\widehat\Dv}_{\max}-k+\tr(\Vv_{\mathcal{A}}\Vv_{\mathcal{A}}^T(\Piv^\star)^T\widehat\Piv)\nonumber\\
    \overset{(b)}{\leq }& |\mathcal{T}|\Bigg(\underbrace{ \frac{4\sqrt{2k}}{\delta}\norm{\Cv_Y-\widehat\Cv_Y}_2}_{\triangleq T_1}-\underbrace{\frac{\norm{\Vv_{\mathcal{A}}^T(\Piv^\star-\widehat\Piv)^T}_F^2}{2\norm{\Piv^\star-\widehat\Piv}_F^2}}_{\triangleq T_2}\Bigg),\nonumber
\end{align}
where $(a)$ follows from H\"older's inequality and \eqref{eq012}, and $(b)$ follows from \eqref{eq_appa12}, \eqref{eq_appa08}, $\norm{\Piv^\star-\widehat\Piv}_F^2=2|\mathcal{T}|$, and $\norm{\Piv^\star-\widehat\Piv}_F^2=2k-2\tr(\Vv_{\mathcal{A}}\Vv_{\mathcal{A}}^T(\Piv^\star)^T\widehat\Piv)$. From the monotonicity property of probability measures, we have
\begin{align}\label{eq_appa06}
	\Pr(\widehat \Piv\neq \Piv^\star)\leq \Pr(T_1\geq T_2).
\end{align}

Let $\epsilon$ be any positive constant. According to \eqref{eq_appa13}, when $t\gtrsim m\ln m$, we have
\begin{align}
    \Pr\left(T_1\geq \frac{8\sqrt{2k}c_0}{\delta}\sqrt{\frac{m(\ln m+u)}{t}}\right)\leq 2e^{-u},\forall u>0.
\end{align}
Setting $c_1=\frac{1}{128c_0^2}$, $t\gtrsim \frac{k\ln m}{\delta^2}$, and $u=\frac{c_1\delta^2t\epsilon^2}{mk}-\ln m> 0$, we have
\begin{align}
      \Pr\left(T_1\geq \epsilon\right)\leq 2me^{-\frac{c_1\delta^2t\epsilon^2}{mk}},\forall \epsilon>0.
\end{align}

Also, setting $[\av_1,\cdots,\av_m]=(\Piv^\star)^T$ and $[\bv_1,\cdots,\bv_m]=\widehat\Piv^T$ in Lemma \ref{lemma_appc}, it  follows that, for any $u>0$:
\begin{align}
    \Pr\left(T_2\leq \frac{k}{2m}(1-u)\right)\leq Ck^{2-C^\prime\frac{ku^2}{\ln(2m)}},
\end{align}
and if $k(1-u)=2m\epsilon$, we have
\begin{align}
    \Pr\left(T_2\leq \epsilon\right)\leq Ck^{2-C^\prime\frac{(k-2m\epsilon)^2}{k\ln(2m)}},\forall \epsilon>0.
\end{align}

Therefore, we have
\begin{align}\label{eq_appc09}
    &\Pr(T_1\geq T_2)=\int \Pr(T_1\geq \epsilon)\Pr(T_2=\epsilon) d\epsilon\nonumber\\
    \leq &\int \Pr(T_1\geq \epsilon)\Pr(T_2\leq\epsilon) d\epsilon\nonumber\\
    \lesssim & mk^2\int e^{-\frac{c_1\delta^2t\epsilon^2}{mk}}k^{-C^\prime\frac{(k-2m\epsilon)^2}{k\ln(2m)}}d\epsilon.
\end{align}

To compute \eqref{eq_appc09}, we need the following Gaussian integral. 
\lemma{Let $a,b,c,$ and $d$ be positive constants with $b\geq 1$, 
\begin{align}
    \int e^{-ax^2}b^{-c(1-dx)^2}dx=\sqrt{\frac{\pi}{a+cd^2\ln b}}e^{-\frac{ac\ln b}{a+cd^2\ln b}}.\label{eq_appc10}
\end{align}
}
\begin{IEEEproof}
    The integral can be derived by employing the Gaussian integral $\int e^{-A(x+B)^2}=\sqrt{\frac{\pi}{A}}$ for $A>0$.
\end{IEEEproof}
Applying \eqref{eq_appc10} to \eqref{eq_appc09} with $a=\frac{c_1\delta^2t}{mk}$, $b=k$, $c=C^\prime k/\ln (2m)$ and $d=\frac{2m}{k}$ and simplifying the expression, we have
\begin{align}\label{eq_appc11}
    &\Pr(\widehat \Piv\neq \Piv^\star)\nonumber\\
    \lesssim& \sqrt{\frac{m^3k^5}{\delta^2 t+m^3\ln k/\ln m}}e^{-c_0\frac{k\delta^2 t\ln k}{\delta^2 t\ln m+m^3\ln k}}\nonumber\\
     \lesssim& \sqrt{\frac{m^3k^5}{m^3/\ln m+\delta^2 t}}e^{-c_0\frac{kt}{\ln m(m^3/\delta^2+t)}}.
\end{align}

\section{Proof of Lemma \ref{lemma1}}\label{app_lemma1}
The result in Lemma \ref{lemma1} directly follows from \cite[Theorem 2]{pmlr-v119-zhang20n}. For completeness, here we show how to adapt the result in \cite[Theorem 2]{pmlr-v119-zhang20n} to obtain Lemma \ref{lemma1}. We note that the notations used in \cite{pmlr-v119-zhang20n} are different from those in this paper. Specifically, the parameters $m$, $n$, and $t$ in this paper are denoted by $n$, $p$, and $m$ in \cite{pmlr-v119-zhang20n}, respectively.

We begin by stating \cite[Theorem 2]{pmlr-v119-zhang20n} using our notations. Suppose the following conditions hold:
\begin{enumerate}
    \item[a.] The entries of $\Av$ are i.i.d. drawn from a standard Gaussian distribution $\Norm(0,1)$. Moreover, each noise vector $\nv_i$ follows the distribution of $\Norm({\bf 0},\sigma^2\Iv_m)$.
    \item[b.] $m\gtrsim n^4(\log n)^4 (\log m)^6$.
    \item[c.] Define $\rho(\Xv^\star)=\norm{\Xv^\star}_F^2/\norm{\Xv^\star}_2^2$. Suppose that $\rho(\Xv^\star)\gtrsim 18/C_0$ and 
    $$\log(\frac{\norm{\Xv^\star}_F^2}{t\sigma^2})\gtrsim\frac{\log(m)}{\rho(\Xv^\star)}+\log \log m$$
    \item[d.] The true permutation $\Piv^\star$ satisfies 
    $$\sum_{i=1}^m \mathbbm{1}_{\{\pi^\star(i)\neq i\}}\leq m/4,$$
    where $\pi^\star(\cdot)$ is the node mapping function corresponding to $\Piv^\star$ and $\mathbbm{1}$ is the indicator function.
\end{enumerate}
Then, according to \cite[Theorem 2]{pmlr-v119-zhang20n}, the permutation $\widehat\Piv_{\text{sp}}$ computed in \eqref{eq17} satisfies
 \begin{align}\label{eq57}
     \Pr(\widehat\Piv_{\text{sp}}\neq \Piv^\star)\lesssim& C_0 e^{-(\log m)^4}+C_1me^{-C_2 t}+C_3me^{-C_4 m}\nonumber\\
     &+C_5e^{-n}+C_6 n^{-2},
 \end{align}
where $C_0,C_1,C_2,C_3,C_4,C_5,C_6$ are some constants.

We show that when $m$, $n$, and $\norm{\Xv^\star}_F^2$ are sufficiently large, Conditions i.-iii. listed in Lemma \ref{lemma1} lead to the above conditions, and the result in \eqref{eq57} yields \eqref{lemma1result}.
\begin{itemize}
    \item Condition a. is the same as Condition i. in Lemma \ref{lemma1}.
    \item When $m$ is sufficiently large such that $\log m>1$, $m\gtrsim n^4(\log n)^4$ suffices $m/(\log m)^6\gtrsim n^4(\log n)^4$. Therefore, Condition b. is satisfied.
    \item Note that $\rho(\Xv^\star)\gtrsim 18/C_0$. Fixing $t$ and $\sigma^2$, Condition c. is satisfied with a sufficiently large $\norm{\Xv^\star}_F^2$.
    \item Note that $\sum_{i=1}^m \mathbbm{1}_{\{\pi^\star(i)\neq i\}}=m-\tr(\Piv^\star)$. Therefore, Condition d. is equivalent to Condition iii. in Lemma \ref{lemma1}.
    \item Finally, we have
    \begin{align*}
        \eqref{eq57}&\lesssim e^{-m}+ m e^{- \min\{C_2,C_4\}\cdot \min\{m,t\}}+C_6 n^{-2}\\
        &\lesssim me^{- c_0^\prime\cdot \min\{m,t\}}+ C_6 n^{-2}.
    \end{align*}
\end{itemize}

Combining all the above results completes the proof.

\section{Proof of Theorem \ref{theorem3}}\label{appc}
Let us denote $\widehat\Zv=\widehat\Piv\Av\widehat\Xv$ and $\Zv^\star=\Piv^\star\Av\Xv^\star$. It follows from \eqref{eq01} that the observation satisfies $\Yv=\Zv^\star+\Nv$ for $\Nv=[\nv_1,\cdots,\nv_t]$. Then, we have
\begin{align}\label{eq_appb01}
 &\norm{\widehat\Zv-\Zv^\star}_F^2=\norm{\widehat\Piv\Av\Av^\dagger\widehat\Piv^T(\Zv^\star+\Nv)-\Zv^\star}_F^2\nonumber\\
 =&\norm{\widehat\Piv\Av\Av^\dagger\widehat\Piv^T\Zv^\star+\widehat\Piv\Av\Av^\dagger\widehat\Piv^T\Nv-\Zv^\star}_F^2\nonumber\\
 \leq & 2\norm{\underbrace{\left(\widehat\Piv\Av\Av^\dagger\widehat\Piv^T\Piv^\star-\Piv^\star\right)\Av\Xv^\star}_{\triangleq \Qv}}_F^2+2\norm{\Av\Av^\dagger\widehat\Piv^T\Nv}_F^2.
\end{align}
We then simplify $\norm{\Qv}_F^2$ as
\begin{align}
    &\norm{\Qv}_F^2=\norm{(\widehat\Piv\Av\Av^\dagger(\widehat\Piv-\Piv^\star)^T\Piv^\star+\widehat\Piv\Av\Av^\dagger-\Piv^\star)\Av\Xv^\star}_F^2\nonumber\\
    \overset{(a)}{=}&\norm{\left(\widehat\Piv\Av\Av^\dagger(\widehat\Piv-\Piv^\star)^T\Piv^\star+\widehat\Piv-\Piv^\star\right)\Av\Xv^\star}_F^2\nonumber\\
    \leq& 2 \norm{(\widehat\Piv-\Piv^\star)\Av\Xv^\star}_F^2+2\norm{\Av\Av^\dagger}_2^2\norm{(\widehat\Piv^T\Piv^\star-\Iv_m)\Av\Xv^\star}_F^2\nonumber\\
   \overset{(b)}{\leq} & 2 \norm{(\widehat\Piv-\Piv^\star)\Av\Xv^\star}_F^2+2\norm{(\widehat\Piv^T\Piv^\star-\Iv_m)\Av\Xv^\star}_F^2\nonumber\\
    =& 4\left(\norm{\Av\Xv^\star}_F^2-\tr(\widehat\Piv^T\Piv^\star\Av\Xv^\star(\Xv^\star)^T\Av^T)\right).\label{eq_appd03}
\end{align}
where $(a)$ follows from $\Av^\dagger\Av=\Iv_n$, $(b)$ is true because $\Av\Av^\dagger$ is a projection matrix with $\norm{\Av\Av^\dagger}_2=1$. 

Denote $\frac{1}{t}\Av\Xv^\star(\Xv^\star)^T\Av^T=\Cv_Z^\star$ and $\Cv_Z=\Av\Cv_X\Av^T$ in \eqref{eq06}. We have $\E[\Cv_Z^\star]=\Cv_Z$. As a result,
\begin{align}\label{eq_appd01}
    \norm{\Av\Xv^\star}_F^2&=t\tr(\Cv_Z)+t\tr(\Cv_Z-\Cv_Z^\star)\nonumber\\
    &\leq t\Bigg(\sum_{i=1}^n\lambda_i+m\norm{\Cv_Z-\Cv_Z^\star}_2\Bigg).
\end{align}
On the other hand,
\begin{align}\label{eq_appd02}
     &\tr(\widehat\Piv^T\Piv^\star\Av\Xv^\star(\Xv^\star)^T\Av^T)\nonumber\\
     \geq& t\tr(\widehat\Piv^T\Piv^\star\Cv_Z)-mt\norm{\Cv_Z-\Cv_Z^\star}_2\nonumber\\
     \geq & t\min_{j\in\mathcal{A}}\lambda_{j}\tr(\widehat\Piv^T\Piv^\star\Vv_{\mathcal{A}}\Vv_{\mathcal{A}}^T)-mt\norm{\Cv_Z-\Cv_Z^\star}_2.
\end{align}
Given $\frac{1}{m}\sum_{i=1}^n\lambda_i\lesssim \frac{n}{m}$ and $k\gtrsim n$ and combining \eqref{eq_appd03}--\eqref{eq_appd02}, it can be verified that
\begin{align*}
    &\frac{\norm{\Qv}_F^2}{mt}\lesssim \frac{n}{m}\left(1-\frac{1}{k}\tr(\widehat\Piv^T\Piv^\star\Vv_{\mathcal{A}}\Vv_{\mathcal{A}}^T)\right)+\norm{\Cv_Z-\Cv_Z^\star}_2.
\end{align*}

Applying Theorem \ref{theorem1} with $\epsilon$ set to $\ln m$, we have
\begin{align}
    \Pr\left(1-\frac{1}{k}\tr(\widehat\Piv^T\Piv^\star\Vv_{\mathcal{A}}\Vv_{\mathcal{A}}^T)\gtrsim \frac{m}{\delta}\sqrt{\frac{ \ln m}{t}}\right)\leq \frac{2}{m}.\nonumber
\end{align}
Similar to \eqref{eq_appa13}, the concentration bound for the sample covariance in \cite[Theorem 5.6.1 and Exercise 5.6.4]{HDimPro} results in
\begin{align}
    \Pr\left(\norm{\Cv_Z-\Cv_Z^\star}_2\gtrsim \sqrt{\frac{\ln m}{t}}\right)\leq \frac{2}{m}.\nonumber
\end{align}
Moreover, since $\Nv$ is a Gaussian matrix, $\Av\Av^\dagger\widehat\Piv^T\Nv$ is still Gaussian distributed with $\E[\norm{\Av\Av^\dagger\widehat\Piv^T\Nv}_F^2]\leq \sigma^2 n/m$.
Applying the tail bound in \cite[Example 2.11]{HDimPro2},
\begin{align}\label{eq_appb02}
    \Pr(\norm{\Av\Av^\dagger\widehat\Piv^T\Nv}_F^2> 2\sigma^2 nt)\leq e^{-\frac{nt}{8}}.
\end{align}

Combining \eqref{eq_appb01}--\eqref{eq_appb02}, with probability at least $1- \frac{4}{m}-e^{-\frac{nt}{8}}$,
\begin{align}
    \frac{\norm{\widehat\Zv-\Zv^\star}_F^2}{mt}\lesssim \frac{n}{\delta}\sqrt{\frac{\ln m}{t}}+\frac{n}{m}\sigma^2.
\end{align}


\bibliographystyle{IEEEtran}
\bibliography{IEEEabrv,ref}

\end{document}